\documentclass[12pt,a4paper]{article}
\usepackage[centertags]{amsmath}
\usepackage{amsfonts,amsthm,amssymb}

\usepackage{amssymb}
\usepackage{amsmath}

\usepackage{graphicx}

\usepackage{setspace}

\def\rit{\hbox{\it I\hskip -2pt  R}}
\newtheorem*{thm1}{Theorem A}
\newtheorem*{thm2}{Theorem B}
\newtheorem*{thm3}{Theorem C}
\newtheorem*{thm0}{Main Theorem 1}
\newtheorem*{thm4}{Theorem D}
\newtheorem*{thm5}{Theorem E}
\newtheorem*{thm6}{Theorem F}
\newtheorem*{thm7}{Theorem G}
\newtheorem*{thm8}{Main Theorem 2}
\newtheorem*{thm9}{Main Theorem 3}
\newcommand{\real}{\mathbb{R}}  

\newcommand{\1}[1]{1_{\{#1\}}}
\newcommand{\In}[2]{#1\!\in\! #2}

\newcommand{\iInN}{\In{i}{N}}

\newcommand{\Norm}[1]{\left\Vert#1\right\Vert}

\newcommand{\cntr}[1]{\begin{center} #1 \end{center}}
\newcommand{\heading}[1]{\cntr{\textbf{#1}}}

\theoremstyle{plain}
\newtheorem{theorem}{Theorem}[section]
\newtheorem{lemma}[theorem]{Lemma}
\newtheorem{claim}[theorem]{Claim}

\theoremstyle{definition}
\newtheorem{definition}{Definition}

\theoremstyle{remark}

\numberwithin{equation}{section}

\begin{document}
\title{Non-Convexity}

\author{Noa Nitzan}
%
%
\maketitle
%
%
%

\section{Introduction}
There are three common measures for evaluating the "non-convexity"
of a set $X \subset \real^d$: \\ $\alpha(X)\,-$ The maximal size of
a visually independent subset of $X$.
\\$\beta(X)\,-$ the minimal size of a collection of seeing subsets of $X$
which covers $X$, or, in other words, the chromatic number of the
invisibility graph of $X$.
\\$\gamma(X)\,-$ the minimum size $k$ such that $X$ can be expressed as the union of $k$ convex sets.

Much effort has been devoted to bounding $\gamma(X)$ in terms of
$\alpha(X)$. In general, there is no such bound, since there exist
planar sets $X$ with $\alpha(X)=3$ but with $\gamma(X)=\infty$, and
there exist closed sets $S \subset \real ^ 4$ with $\alpha(S)=2$ and
$\gamma(S)=\infty$ (even $\beta(S)=\infty$). In the specific case of
closed sets in the plane, the situation is different.

Valentine [1957] proved that for closed $S\subset\real^2$,
$\alpha(S)=2$ implies $\gamma(S)\leq 3$. Eggleston [1974] proved
that for compact $S\subset\real^2$, $\alpha(S)<\infty$ implies
$\gamma(S)<\infty$. Breen and Kay [1976] were the first to find an
upper bound for $\gamma$ in terms of $\alpha$. They proved that for
closed $S\subset\real^2$, if $\alpha(S)=m$ then $\gamma(S)\leq
m^3\cdot2^m$. Later on, Perles and Shelah [1990]  improved this
upper bound to $m^6$, and Matou$\check{\text{s}}$ek and Valtr [1999]
obtained the best upper bound know today, $18m^3$. In the same
paper, M. and V. give examples of closed planar sets $S$ with
$\alpha(S)=m$ and $\gamma(S)=cm^2$

There has also been some success in bounding $\gamma(X)$ for certain
cases of $X$ not necessarily closed. Breen [1974] claims that for
$X\subset \real^2$, $\alpha(X)=2$ implies $\gamma(X) \leq 6$.
Another result is of Matou$\check{\text{s}}$ek and Valtr [1999] who
proved that for $X \subset \real^2$ with finite $\alpha(X)=m$, if
$X$ is starshaped then $\gamma(X)\leq 2m^2$. They also proved that
for $X \subset \real^2$, if $\real^2 \setminus X$ has no isolated
points then $\gamma(X)\leq m^4$.

In this work we shall focus on the case of $X \subset \real^2$ with
$\alpha(X)=2$. We wish to complete the work of Breen [1974], and
give a detailed proof of the theorem claimed by Breen
($\alpha(X)=2\Rightarrow \gamma(X)\leq6)$. We intend to determine
the maximum possible value of $\gamma(X)$ (assuming $X \subset
\real^2$ and $\alpha(X)=2$) under a variety of side conditions,
pertaining to the location (within $\textrm{cl}X$) of the points of
$\textrm{cl}X$ that are missing in $X$. We produce examples for all
cases under discussion, showing that the bounds obtained are tight.

\section{Definitions and Notations}

Given $X\subseteq\real^2$, we say that two points $u,v\in \real^2$
$\mathbf{see}$ $\mathbf{each}$ $\mathbf{other}$ via $X$ if the open
interval $(u,v)$ is included in $X$. (This applies even if the
points $u,v$ are not in $X$)
\\$A$ is a $\mathbf{seeing}$ $\mathbf{subset}$ of $X$ if $A \subset X$ and every two points of $A$ see each
other via $X$.
\\A subset of $X$ is $\mathbf{visually}$ $\mathbf{independent}$ if no two of its
points see each other via $X$.
\\Define the $\mathbf{invisibility}$
$\mathbf{graph}$ of $X$ as the graph $G(X)$ with vertex set $X$ and
with $u,v\in X$ connected by an edge iff $[u,v] \nsubseteq X$.

We now define the 3 most common ``measures of non-convexity'' of
$X$: (These are the notations found in the literature which we
prefer.)

$\alpha(X)\,-$ The supremum of cardinalities of all visually
independent subsets of $X$. That is, the clique number of the graph
$G(X)$.

$\beta(X)\,-$ The chromatic number of $G(X)$. In other words, the smallest \\cardinality of a collection of seeing subsets of $X$
that covers $X$.

$\gamma(X)\,-$ the smallest cardinality $k$ such that $X$ can be expressed as the union of $k$ convex sets.
\\It is easy to see that $\alpha(X)\leq \beta(X)\leq \gamma(X)$.

The following notations will be used throughout this paper: For $X
\subset \real^2$, define $S=\textrm{cl}X$. We shall write
$M=S\setminus X$ $\;$($M$ is the set of points of $S$ missing in
$X$). We split $M$ into two parts $M=M_b \cup M_i$, where $M_i =M
\cap \textrm{int} S$ and $M_b=M \cap \textrm{bd}S$.

$S$ is $\mathbf{locally}$ $\mathbf{convex}$ at a point $x$ if $x \in
S$ and $x$ has a neighborhood $U$ such that $S \cap U$ is convex. We
denote by $Q$ ($=\textrm{lnc}S$ the set of points of local
non-convexity ($\textrm{lnc}$ points) of $S$. These are the points
where $S$ fails to be locally convex.

We say that $S$ is $2$-dimensional at a point $p$ if $p \in
\textrm{cl}(\textrm{int}S)$.

$A \subset \real^d$ is an $L_2$-set if every two points of $A$ can
be connected by a polygonal line of at most 2 edges within $A$.

Given a subset $S_0 \subset S$, we say that $S_0$ is convex relative to $S$ if for every $x,y \in S_0$, if $[x,y] \subset S$ then $[x,y] \subset S_0$.
\section{Results}
Throughout the following theorems we assume that $X$ is a planar
set, $\alpha(X)\leq 2$ and that $S=\textrm{cl}X$.
\begin{thm0}
$max\{\gamma(X):X \subset \real^2,  \alpha(X)\leq2 \}=6$
\end{thm0}

\begin{thm1}
If $X$ is  not an $L_2$-set (in particular, if $X$ is not
connected), then $\gamma(X)=2$. (In this theorem, $\real^2$ can be
replaced by an arbitrary real vector space)
\end{thm1}

\begin{thm2}
If $S$ is  not $2$-dimensional at some point, then $\gamma(X)\leq
2$.
\end{thm2}

\begin{thm3}
If $\;|M_i|>1$, then $\gamma(X)\leq 3$. The number three is best
possible, even when $S$ is convex. If, in addition, $M_b=\phi$ or
$M_b=$\textrm{bd}$S$ then $\gamma(X)=2$.
\end{thm3}

\begin{thm4}
If $\;|M_i|=1$ and $M_b=\phi$ or $M_b=\textrm{bd}S$, then $\gamma(X)\leq 4$. The number four is best possible.
\end{thm4}

\begin{thm5}
If $\;M_i=\phi$ then $\gamma(X)\leq 3$. The number three is best
possible, even when $S$ is convex.
\end{thm5}

\begin{thm6}
If $\;|M_i|=1$ then $\gamma(X)\leq 6$. The
number six is best possible.
\end{thm6}

\begin{thm7}
If $\;|M_i|=1$ and $S$ is convex then $\gamma(X)\leq 4$. The number
four is best possible. If, in addition, $M_b=\phi$ or
$M_b=\textrm{bd}S$ then $\gamma(X)=2$.
\end{thm7}
\begin{thm8}
$max\{\beta(X):X \subset \real^2,  \alpha(X)\leq2 \}=4$
\end{thm8}

\begin{thm9}
$max\{\gamma(X):X \subset \real^2,  \beta(X)=2 \}=4$
\end{thm9}

Table 1 summarizes all the cases above. In each box appears $max\,
\gamma(X)$ under the conditions of that box. The number in
parentheses is $max \,\gamma(X)$ under the conditions of the box
together with the extra assumption that $S$ is convex.

Much of the material contained in this paper can be summarized in
the following extension of Valentine's Theorem:

\begin{theorem}
If $X \subset \real^2$, $\alpha(X) \leq2$, and the complement
$\real^2 \setminus X$ has no one-pointed components, then $\gamma(X)
\leq 3$.
\end{theorem}
\section{Proof of Theorem A}
As $X$ is  not an $L_2$-set, there are two points $a,b \in X$ which
cannot be connected by a polygonal line of less than 3 edges within
$X$. In other words, there is no point in $X$ that sees both $a$ and
$b$. Define $A=st(a)=\{x \in X:[a,x]\subset X\}$, $B=st(b)$. Notice
that the sets $A$ and $B$ are disjoint. We show now that $A$ is
convex: Take $p,q \in A$ where $p=a+u$, $q=a+v$. For every $0<
\theta \leq 1$, $[a+\theta u, a+\theta v ]\subset X$ because
otherwise $\{ a+\theta u, a+\theta v, b \}$ is a visually
independent set, a contradiction. Hence the full triangle $[a,p,q]$
is included in $X$, so $a$ sees via $X$ every point in $[p,q]$,
which means that $[p,q] \subset A$, so $A$ is convex. Similarly, $B$
is convex. Now, for every $x \in X$, $x \in A \cup B$, because
otherwise, $\{a,b,x\}$ is a visually independent set, a
contradiction.
We can now conclude that $X$ is the union of two disjoint convex sets.

\section{Proof of Theorem B}

$S$ is a closed set in the plane and therefore, according to
Valentine [1957], is a union of at most three convex sets:
$S=\cup_{i=1}^n{C_i}$ where $1\leq n \leq 3$. As $S$ is closed, we
can assume that for each $i$, $C_i$ is closed. In addition, we will
assume that none of these convex sets is included in the union of
the others.

If each $C_i$ is $2$-dimensional then $S$ is $2$-dimensional. Assume therefore, w.l.o.g., $C_1$ is not $2$-dimensional. If $C_1$ is of
dimension $0$, then $S$ is not connected and therefore $X$ is not
connected, so by theorem $A$, $X$ is the union of 2 convex sets.
Otherwise, $C_1$ is of dimension $1$, meaning that $C_1$ is a
segment. There is point $p \in C_1$ such that $p \notin C_2\cup
C_3$. $C_2,C_3$ are closed so the distances $d(p,C_2),d(p,C_3)$ are
positive. Hence there is a neighborhood $\,U$ of $p$ such that
$U\cap S$ is a segment. Define $L$ to be the line containing this
segment. Denote by $L_+$, $L_-$ the open half-planes determined by
$L$.

Define $C=\textrm{conv}(X\setminus L)$. We wish to show that $C \subset X$:
Every two points in $X\cap L_+$ do not see $p$ via $X$, therefore, since $\alpha(X)=2$, they see each other via $X \cap L_+$. Hence, $X\cap
L_+$ is convex. By the same argument, $X\cap L_-$ is convex, so
$C=\textrm{conv}(X\setminus L)=\textrm{conv}((X\cap L_+)\cup (X\cap
L_-))=\cup \{[a,b]:{a\in {X \cap L_+}, b \in {X \cap L_-}}\}$. The
point $p$ does not see any $a\in L_+ \cap X$ or $b \in L_-\cap X$,
therefore, again, as $\alpha(X)=2$, for any such $a,b$,
$[a,b]\subset X$.
 This implies that $C \subset X$.

It remains to deal with the set $L\cap X$. Since $\alpha(X)=2$ and
$L$ is convex, $\alpha(X \cap L)\leq 2$. If $X \cap L$ is convex, we
are done. Otherwise, $X \cap L$ is the disjoint union of two
nonempty convex sets $A$, $B$, where say, $p \in A$. If $C=\phi$
then we are done, so assume $C\neq \phi$.

In order to complete the proof, we would like to show that $\textrm{conv}(B
\cup C)\subset X$. Since both $B$ and $C$ are convex,
$\textrm{conv}(B \cup C)=\cup\{ [b,c]:{b\in B, c \in C}\}$.
Let $b \in B$ and $c \in C$:
\\Case 1: If $c\notin L$ then $[p,c]\nsubseteq X$ and $[p,b] \nsubseteq X$, hence
$[b,c] \subset X$.

Case 2: If $c \in L$ then $c \in [a_+, a_-]$, where $a_+ \in X \cap L_+$ and
$a_- \in X \cap L_-$. The points $a_+, a_-, b$ do not see $p$ via
$X$, therefore, and since $\alpha(X)=2$, $[a_+,b] \subset X$ and
$[a_-,b] \subset X$. Now, each point in $[c,b)$ is in the convex
hull of a point in $[a_+,b)$ and a point in $[a_-,b)$ and therefore
is in $X$ (again, these two points do not see $p$, and
$\alpha(x)=2$).

This establishes that $\textrm{conv}(B \cup C)\subset
X$, which implies that $X$ is the union of two convex sets: $A$, the
component of $p$ in $L$, and $\textrm{conv}(B \cup C)$(=$B \cup C$).

\section{Proof of Theorem C}

Coming to prove theorem C, we shall first show that $M_i$ contains a
segment. Suppose $x,y \in M_i$, $x \neq y$, and let  $L$ be the line
spanned by $ x,y$. As $x,y \in \textrm{int}S$, both have
neighborhoods $U_x,U_y$ in $S$. The intersection of $L\setminus
\{x,y\}$ with these two neighborhoods provides 4 segments. These
segments lie in the three components of $L \setminus \{x,y\}$, and
therefore at least one of them is disjoint from $X$. Therefore $M_i$
contains a segment, call it $I$.

Denote by $L_+$, $L_-$ the open half-planes determined by $L$. Define $X_+=X \cap L_+, X_-=X
\cap L_-$ and $S_+=\textrm{cl}(X_+), S_-=\textrm{cl}(X_-)$. Next we
show that  $X_+$ is convex: Take $p,q \in X_+$. There is $y \in X_-
$, close enough to the center of $I$, such that both segments
$(p,y), (q,y)$ intersect $I$, meaning that $y$ sees neither $p$ nor
$q$ via $X$, and therefore $[p,q] \subset X$, hence $[p,q] \subset
X_+$. Similarly, $X_-$ is convex.

We wish to show now that $M_i \subset L$: Note that $S \cap L_+
\subset \textrm{cl}X_+$, and therefore $L_+ \cap
\textrm{int}S=\textrm{int}(S \cap L_+) \subset
\textrm{int}\,\textrm{cl}X_+=\textrm{int}X_+ \subset X$. Therefore
$L_+ \cap M_i= \phi$. Similarly, $L_- \cap M_i= \phi$, hence $M_i
\subset L$.

Define: $B_+=S_+ \cap L$, $B_-=S_- \cap L$. $B_+$ is an edge of $S_+$
(think of it as the base of $S_+$). If a point $u$ lies in
$\textrm{rel}\,\textrm{int}B_+$, then $u$ sees every point of $S
\cap L_+$ via $\textrm{int}S_+$, hence via $X$. Thus, a point $u \in
X \cap B_+$ may fail to see some point of $X_+$ via $X$ only if $u$
is an endpoint of $B_+$. Similarly for $B_-$ and $X_-$.

If $X$ contains a point $y$  that is in $L\setminus(B_- \cup B_+)$,
then $S$ is not $2$-dimensional at $y$, and therefore $\gamma(X)=2$,
by Theorem B. Assume therefore that $X \cap L \subset B_- \cup B_+$.
Note that the segment $I$($\subset M_i$) lies in $B_- \cap B_+$.
\\Next we show that $\gamma(X)=2$, unless $X \cap L \subset B_- \cap
B_+$. Assume $X \cap L \nsubseteq B_- \cap B_+$. Pick a point $y \in
X \cap L \setminus (B_- \cap B_+)$. Think of $L$ as a horizontal
line, and suppose, w.l.o.g., that $y \notin B_-$, and that $y$ is to
the right of $B_-$.

Denote by $L_2$ the component of $y$ in $X \cap L$. $L_1$ is the
other component of $X \cap L$, if $X \cap L$ is not convex. If $X
\cap L$ is convex, then $L_1=\phi$.

Clearly, $y$ does not see any point of $L_1$ via $X$. Since $y
\notin B_-$, $y$ doesn't see any point of $X_-$ via $X$. Since
$\alpha(X)=2$, every point of $L_1$ sees every point of $X_-$ via
$X$, hence via $X_-$. In other words, $L_1 \cup X_-$ is convex. But,
as we shall see immediately, $L_2 \cup X_+$ is also convex. Indeed,
consider a point $x \in X_+$ and a point $y' \in L_2$, to the right
of $y$ ($y'=y$ included). $x$ doesn't see via $X$ some point $z \in
X_-$, that lies beyond $I$. $y' \notin  B_-$ and therefore doesn't
see via $X$ any point in $X_-$. It follows that $y'$ sees $x$ via
$X$. Now consider a point $y'' \in L_2$, strictly to the left of
$y$. Since $y'' \in L_2$ lies to the right of $I$, and $I \subset
B_+$, we conclude that $y'' \in \textrm{rel}\,\textrm{int}B_+$, and
therefore sees via $X_+$ every point of $X_+$. We can now represent
$X$ as the union of two convex sets: $X= (L_1 \cup X_-) \cup (L_2
\cup X_+)$.

Assume from now on that $X \cap L \subset B_- \cap B_+$. Let us first dispose of the case where $M_b=\phi$ or
$M_b=\textrm{bd}S$.

$M_b=\phi$: If $c \in X \cap L$$(\subset B_- \cap B_+)$ then $c$
sees every point of $X_+$ via $S \cap L_+$, which is a subset of
$X_+$. Same for $X_-$ and $L_-$. Denote by $L_1,L_2$ the components of $X \cap L$. ($L_1= \phi$ if
$X \cap L$ is convex). Then $X$ is the union of the two convex sets
$L_1 \cup X_+$, $L_2 \cup X_-$.

$M_b=\textrm{bd}S$: If $x \in X \cap L$ then $x \in
\textrm{rel}\,\textrm{int}B_+$. (The endpoints of $B_+$ are boundary
points of $S$ and therefore not in $X$.) Similarly, $x \in
\textrm{rel}\,\textrm{int}B_-$. Define $L_1,L_2$ as above. Then $X$ is again the union of the two
convex sets $L_1 \cup X_+$ and $L_2 \cup X_-$.

Now we return to the general case: $|M_i|>1$ and $M_b$
unrestricted, and try to show that $\gamma(X)\leq 3$.

If $X \cap L$ is convex then $X$ is
the union of three convex sets and we are done. Assume $X \cap L$ is
not convex, so it is composed of two non-empty components $L_1,L_2$,
where $L_1$ is to the left of $L_2$.
If $L_1$ has no left endpoint, then $L_1 \subset
\textrm{rel\,int}B_+$, and therefore $X$ is the union of the three
convex sets $\,L_1 \cup X_+,X_-,L_2$. The same argument works when
$L_1$ has a left endpoint $c_1$, but $c_1$ is not the left endpoint
of $B_+$. We can repeat this argument with $B_-,X_-$ instead of
$B_+,X_+$, and also with $L_2$ instead of $L_1$.

Assume, therefore, that $L_1$ has a left endpoint $c_1$, $L_2$ has
a right endpoint $c_2$, and $B_+=B_-=[c_1,c_2]$. The point $c_1$
still sees  every point of $S \cap L_+$ via $\textrm{int}S_+(\subset
X)$ unless $S_+$ has an edge $C_{1_+}$ with endpoint $c_1$, other
than $B_+$. Assume, therefore that $S_+$ has such an edge $C_{1_+}$,
and by the same token, that $S_+$ has has an edge $C_{2_+}$ with
endpoint $c_2$, other than $B_+$ (see Figure 1).
If $X \cap C_{1_+}$ is convex then $c_1$ still sees every point of
$X_+$ via $X_+$, and thus $L_1 \cup X_+$ is again convex, as before.

Assume therefore that $X \cap C_{1_+}$is not convex. It is the
union of $\{c_1\}(=C_{1_+} \cap L)$ and the convex set $C_{1_+} \cap
X_+$. By the same token, assume that $X \cap C_{2_+}$ is not convex.
It follows that $X \cap C_{1_+} \cap C_{2_+} = \phi$ since a point
$z \in X \cap C_{1_+} \cap C_{2_+}$ would form a 3-circuit of
invisibility with $c_1$ and $c_2$.

We could play the same game with $X_-$, but this is not necessary,
since $X$ is the union of the three convex sets $X_-$,
$(X_+\setminus C_{1_+}) \cup L_1$, $(X_+ \setminus C_{2_+}) \cup
L_2$.

Examples 1,2 show that the number three is best possible: We describe two sets $X_1,X_2 \subset \real^2$ with $|M_i|>1$ and show that $\alpha$ of each set is 2 and that $\gamma$ of each set is not less than 3.
Notice that $X_1 \cap L$ is not convex, while $X_2 \cap L$ is convex.
\\\underline{\textbf{Example 1}}:

Let $P$ be a regular hexagon with center $O$ and vertices $p_0,p_1,...,p_5$. Take $[a,b]$ to be a short
 segment lying on $[p_5,p_2]$ with $O$ in its center. We define
$X_1=P \setminus ((p_5, p_0 ) \cup (p_1,p_2) \cup [a,b])$ (see
Figure 2).

$\alpha(X_1)=2$: The set $X_1 \setminus \{p_0\}$ is the union of two convex sets. The same holds for $X_1 \setminus \{p_1\}$. Therefore, if there
is a 3-circuit of invisibility in $X_1$, it must contain both $p_0$
and $p_1$. But these two points see each other via $X_1$.

$\gamma(X_1)\geq 3$ since, as shown in Figure 2, there is a 5-circuit of invisibility.
\\\underline{\textbf{Example 2}}:

Let $P$ be as above and take $[c,d]$ to be a short segment lying on $[p_5,p_2]$ with $c=p_5$. Define $X_2=P \setminus ((p_1, p_2 ) \cup (p_2,p_3) \cup [c,d])$ (see Figure 2).

$\alpha(X_2)=2$: The set $X_2 \setminus \{p_2\}$ is the union of two convex sets.
 Therefore, if there is a 3-circuit of invisibility in $X_2$, it must contain $p_2$.  But there are only two points which $p_2$ doesn't see via $X_2$: $p_1$ and $p_3$, and fortunately $[p_1,p_3] \subset X_2$.

 $\gamma(X_2)\geq 3$ since, as shown in Figure 2, there is a 5-circuit of invisibility.

\section{Proof of Theorem D}

\begin{lemma}\label{noa3}
$M_i \subset \textrm{ker}S$
\end{lemma}
\begin{proof}
Assume $x \in M_i$ and suppose there is a point $y \in S$ such that
$[x,y] \nsubseteq S$. In other words, there is a point $z \in [x,y]$
such that $z \notin S$. As $S$ is closed, there is a neighborhood
$U$ of $z$, disjoint from $S$. $y \in S=\textrm{cl}X$, so there is a
point $y' \in X$, close to $y$, satisfying $[x,y']\cap U \neq \phi$.
$x \in \textrm{int}S$, so there is a neighborhood $V\subset
\textrm{int}S$ of $x$ such that all points in $V$ do not see $y'$
via $S$. It is possible to choose three points $a,b,c$ in $V$ such
that $x \in \textrm{int\,conv}(a,b,c)$. A slight perturbation of
these three points will lead to $a',b',c' \in X\cap V$, such that
all three points do not see $y'$ via $X$ and $x \in \textrm{int
\,conv}(a',b',c')$. Now, since $y'$ doesn't see any of $a',b',c'$
via $X$, the segments $[a',b'],[b',c'],[c',a']$ are in $X$ (this is
true since $\alpha(X)= 2$). Therefore, it is possible to find a
3-circuit of invisibility in $X$: $y'$, and the two points of
intersection of any line through $x$ with the boundary of the
triangle $[a',b',c']$. This contradicts $\alpha(X)= 2$.
\end{proof}
Before we prove theorem $D$ we shall quote a result
of Breen and Kay [1976]: Let $S\subset \real^2$ be a closed set with
finite $\alpha(S)$. If $S$ is starshaped with respect to a point
that lies on a line that supports $S$, then $\gamma(S)=\alpha(S)$.

We shall now return to the proof of theorem $D$.

\underline{If $M_b= \phi$}: Assume
$M_i=\{ (0,0)\}$, thus $X=S \setminus \{(0,0)\}$. Define $S_+=S \cap
\{(x,y)|y\geq 0\}$ and $S_-=S \cap \{(x,y)|y\leq 0\}$. $S=S_+ \cup
S_-$. Since $S_+$ is the intersection of $S$ with a convex set,
$\alpha(S_+)\leq 2$. According to the lemma above, $S_+$ is
starshaped with respect to $(0,0)$ and therefore, due to Breen and
Kay, $\gamma(S_+)=\alpha(S_+)\leq 2$, so $S_+=A \cup B$ for convex
sets $A,B$. Similarly, $\gamma(S_-)\leq 2$, so $S_-=C\cup D$ for
convex sets $C,D$, hence $S=A \cup B \cup C  \cup D$.

Define $\,T_+=\{(x,y)|y>0 \vee (y=0 \wedge x>0) \}$
and $\,T_-=\{(x,y)|y<0 \vee (y=0 \wedge x<0) \}$. $T_+$ and $T_-$
are both convex, and $\;T_+ \cup T_- = \real^2 \setminus \{(0,0)\}$.

We wish to show that in this case, when $M_b= \phi$, $X$ is the following union of four convex sets: $X=(A
\cap T_+) \cup (B \cap T_+) \cup (C\cap T_-) \cup (D \cap T_-)$.
It is clear that the union of the four sets is included in $X$. We shall see the opposite inclusion:
Suppose $p=(x,y) \in X$. Note that $p \neq (0,0)$. \\If $y>0$, or if
$y=0$ and $x>0$, then $p\in S_+ \cap T_+=(A \cap T_+) \cup (B \cap
T_+)$.
\\If $y<0$, or if
$y=0$ and $x<0$, then $p\in S_- \cap T_-=(C \cap T_-) \cup (D \cap
T_-)$.

There is an alternative proof of Theorem $D$ which avoids the result
of Breen and Kay quoted above. Instead, it uses the necessary
condition for $\gamma(S)=3$ in Valentine's Theorem. We shall use
this alternative approach in the proof of the slightly more
complicated case $M_b=\textrm{bd}S$.

\underline{If $M_b=\textrm{bd}S$}: Assume $M_i=\{p\}$, $p \in
\textrm{int}S$. Then $X=S\setminus \textrm{bd}S\setminus
\{p\}=\textrm{int}S \setminus \{p\}$. Assume $\gamma(S)=k$, $1\leq
k\leq 3$. Then $S$ is the union of $k$ closed convex sets $C_i$
$(i=1,...,k)$. Replacing $C_i$ by $\textrm{conv}(\textrm{ker}S \cup
C_i)$, if necessary, we may (and shall) assume that $\textrm{ker}S
\subset C_i$ $(i=1,...,k)$, and clearly, $\cap_{i=1}^k
C_i=\textrm{ker}S$. Now, consider the following cases:

1) $k=1$, i.e., $S$ is convex. So is $\textrm{int}S$, and
$\textrm{int}S\setminus \{p\}$ is the union of two convex sets.

2) $k=2$ and $\textrm{dim\,ker}S=2$. From $\textrm{int}C_1 \cap
\textrm{int}C_2\neq \phi\,$ it follows that
$\textrm{int}S=\textrm{int}C_1 \cup \textrm{int}C_2$. (Clearly,
$\textrm{int}C_1\cup \textrm{int}C_2 \subset \textrm{int}S$.
Conversely, if $x \in \textrm{int}S$, and $z \in \textrm{int}C_1
\cap \textrm{int}C_2$, then, for some sufficiently small
$\epsilon>0$, $x'=(1+\epsilon)x- \epsilon z \in S =C_1 \cup C_2$.
Suppose, say, that $x' \in C_1$, then $x \in (x',z] \subset
\textrm{int}C_1$.)

Therefore, $X=\textrm{int}S \setminus \{p\}=(\textrm{int}C_1
\setminus \{p\}) \cup (\textrm{int}C_2 \setminus \{p\})$ is the
union of at most 4 convex sets. Example 3 shows that sometimes $X$
is not the union of fewer than 4 convex sets.

3) $k=2$ and $\textrm{dim\,ker}S=1$. Put $K=\textrm{ker}S$,
$L=\textrm{aff}K$, and let $L_+$, $L_-$ be the two closed
half-planes bounded by $L$. Then $K$ is a closed line segment (or a
ray) within $L$.

Since $S$ has no lnc points outside $L$, it follows that the sets
$S_+=S \cap L_+$, $S_-=S \cap L_-$ are convex, $S=S_+ \cup S_-$ and
$S_+ \cap S_-=K$. It follows easily that
$\textrm{int}S=\textrm{int}S_+ \cup \textrm{int}S_- \cup
\textrm{rel\,int}K$.

Now $M_i=\{p\}$, where $p \in \textrm{int}S \cap \textrm{ker}S$.
Thus $p \in \textrm{rel\,int}K$. The point $p$ divides
$\textrm{rel\,int}K$ into two (one-dimensional) convex sets $K_1,K_2
\subset L$, and $X=\textrm{int}S\setminus \{p\}$ is the union of two
convex sets: $X=(K_1 \cup \textrm{int}S_-)\cup (K_2 \cup
\textrm{int}S_+)$.

4) $k=2$ and $\textrm{dim\,ker}S\leq 0$. Thus $S=C_1 \cup C_2$,
where $C_1$ and $C_2$ are closed convex sets,
$\textrm{dim}C_1=\textrm{dim}C_2=2$ (otherwise, according to Theorem
B, $\gamma(X) \leq 2$). $|\textrm{ker}S|=|C_1 \cap C_2| \leq 1$. If
$C_1 \cap C_2 = \phi$ then there is no room for $p \in M_i \subset
\textrm{ker}S$. If $|C_1 \cap C_2|=1$ then $C_1 \cap C_2= \{p\}$. In
this case, $X=\textrm{int}C_1 \cup \textrm{int}C_2$ is the
(disconnected) disjoint union of two convex sets.

5) $k=3$. From $\alpha(S)=2$ and $\gamma(S)=3$ it follows (due to
Valentine[57]) that $S$ is the union of an odd-sided convex polygon
$P=\textrm{conv}Q$ (where $Q=\textrm{lnc}S=\{q_1,...,q_m\}$, $m\geq
3, m$ odd) and $m$ leaves $W_1,...,W_m$. Each leaf $W_i$ is a closed
convex set that includes the edge $[q_i,q_{i+1}]$ of $Q$ (where
$q_{m+1}=q_1$), and that lies beyond that edge and beneath all other
edges of $P$. Note the subset $P \cup \cup_{i=1}^{m-1} W_i$ of $S$
is the union of two convex sets: $P \cup \{W_i:1 \leq i <m, i \text{
odd} \}$ and $P \cup \{W_i:i \text{ even} \}$.

The missing point $p \in \textrm{int}S \cap \textrm{ker}S$ may lie
in $P$ or outside $P$. It is certainly not a vertex of $P$. Pass a
line $L$ through $p$ that passes through $\textrm{int}P$ but misses
all vertices of $P$. Denote by $L_+,L_-$ the two closed half-planes
determined by $L$, and define:

$S_+=S \cap L_+$     $\;\;\;\;S_-=S \cap L_-$

$P_+=P \cap L_+$     $\;\;\;\;P_-=P \cap L_-$

$W_{i_+}=W_i \cap L_+$     $\;\;\;\;W_{i_-}=W_i \cap L_-$
$(i=1,2,...,m)$

$P_+$ is a convex polygon, even-sided or odd-sided. For $i=1,...,m$,
$W_{i_+}$ is either empty, or a closed, convex leaf that sits on an
edge of $P_+$ and lies beneath all other edges of $P_+$. But there
is no leaf sitting on the edge $P \cap L$ of $P_+$ (see Figure 3).
 It follows that
$S_+$ is the union of two closed convex sets: $S_+=C_{1_+} \cup
C_{2_+}$. These two convex sets can be extended to include the
convex kernel of $S_+$. We shall therefore assume that $\{p\} \cup
P_+ \subset C_{i_+}$ for $i=1,2$. The same argument, with $+$
replaced by $-$, applies to $S_-$.

Denote by $K_1,K_2$ the two components of the set $S \cap L
\setminus \{p\}$. Then $\textrm{int}(S \setminus
\{p\})=\textrm{int}S_+ \cup \textrm{int}S_- \cup \textrm{relint}K_1
 \cup \textrm{relint}K_2$.

Since $K_1 \subset S_+=C_{1_+} \cup C_{2_+}$, and $p \in C_{1_+}
\cap C_{2_+}$, one of the sets $C_{1_+}, C_{2_+}$, say $C_{1_+}$,
must include $K_1$. By the same argument, applied to $K_2$ and
$S_-$, one of the sets $C_{1_-},C_{2_-}$, say $C_{2_-}$, must
include $K_2$. Thus $X=\textrm{int}S\setminus \{p\}$ is the union of
the four convex sets:  $\textrm{int}C_{1_+} \cup
\textrm{relint}K_1,\;\textrm{int}C_{2_+},\;\textrm{int}C_{1_-},$
$\textrm{int}C_{2_-} \cup \textrm{relint}K_2$.

A well known example shows that the number four is best possible.
\\\underline{\textbf{Example 3}}:

We describe a set $X \subset \real^2$ with $|M_i|=1$ and $M_b = \phi$, and show that $\alpha(X)=2$ and $\gamma(X)\geq 4$.
Let $X$ be a closed Star of David with its center $O$ removed.
Denote by $p_0,p_1,...,p_5$ the vertices of the Star of David which
are locally convex (see the left side of Figure 4).

$\alpha(X)=2$: The right side of Figure 4 shows a representation of $X$ as the union of two seeing subset, hence
$\beta(X)=2$. Therefore $\alpha(X)=2$.

$\gamma(X)\geq4$: Let $C$ be a convex subset of $X$. Define $A=\{p_0,p_2,p_4\}$ and $B=\{p_1,p_3,p_5\}$. No point of $A$
sees any point of $B$, therefore $C$ cannot contain points of both
sets. $A \nsubseteq C$ and $B \nsubseteq C$, since $O \in C$.
Therefore, $C$ contains at most two points of $A$ or two points of
$B$. Hence, for each of the sets $A,B$, at least two convex sets are
needed in order to cover it. It follows that $\gamma(X)\geq 4$. One
can easily represent $X$ as a disjoint union of 4 convex sets.

For an example with $|M_i|=1$ and $M_b = \textrm{bd}S$, take the same set $X$ and  remove its boundary.
\section{Proof of Theorem E}
Denote $K=\textrm{ker}S$. We prove theorem E by classifying
$\gamma(X)$ in terms of the dimension of $K$. We first need the
following lemmata:
\begin{lemma}\label{noa4} $S$ has no triangular holes.
\end{lemma}
\begin{proof}
If $S$ is convex, then of course, there are no holes in $S$.

If $S$ is not connected, Then $S$ is the union of two disjoint, convex
sets, and again there are no holes in $S$.

Otherwise, if $S$ is connected but not convex, then according to
Tietze's Theorem, $S$ contains an lnc point $q$. According to
Valentine [1957], $q \in \textrm{ker}S$. In other words, $S$ is
starshaped with respect to $q$. But a starshaped set has no holes.
\end{proof}
\begin{lemma}\label{noa1} If $M_i
= \phi$, then $\beta(X)=\gamma(X)$.
\end{lemma}
\begin{proof}
It suffices to show that if $A$ is a seeing subset of $X$ then
$\textrm{conv}A \subset X$. Every point in $\textrm{conv}A$ is a
convex combination of at most 3 points of $A$. If $x$ is a convex
combination of 2 points of $A$, then $x \in X$. Assume $x$ is a
convex combination of 3 affinely independent points $a,b,c \in A$.
The edges of the triangle $\Delta=[a,b,c]$ are included in $X$.
According to Lemma \ref{noa4}, $S$ has no triangular holes,
therefore $\Delta \subset S$. This implies $x \in \textrm{int}\Delta
\subset \textrm{int}S \subset X$. (The last inclusion is due to
$M_i= \phi$.)
\end{proof}
In view of Lemma \ref{noa1}, we only have to find how many seeing
subsets of $X$ are needed in order to cover $X$.
\\

Case 1:
$K=\phi$. If the kernel is empty then $S$ is not connected and so is
$X$, so by theorem A, $X$ is the union of two convex sets.
\\

Case 2: $\textrm{dim}K=0$. We will show that in this case, $X$ is
the union of two convex sets. When $|K|=1$, then according to the
 the proof in Valentine [1957], $S$ is the union of two convex sets. We can assume
that $S$ is the union of two closed, convex sets $A,B$, both
containing $K=\{q\}$, so $A \cap B =\{q\}$. We can also assume that
both $A,B$ are of full dimension, otherwise, we are back to the case
of Theorem $B$. If $q \notin X$ then $X$ is not connected, so assume
$q \in X$.

We claim that $X \cap (A\setminus \{q\})$ is a seeing subset of $X$.
Indeed, suppose $x,y \in X \cap (A\setminus \{q\})$. Note that if $b
\in B \cap X$, then $x$ cannot see $b$ via $X$ (even via $S$) unless
$q \in [x,b]$. Similarly for $y$. Chooses a point $b \in B$ that is
not collinear with $x,q$, nor with $y,q$ ($\textrm{dim}B=2$). Then
$b$ sees neither $x$ nor $y$ via $X$, and therefore $[x,y] \subset
X$. By the same token, $X \cap (B \setminus \{q\})$ is also a seeing
subset  of $X$.

We still have to take care of the point $q$. We would like to add
$q$ to either $X \cap (A \setminus \{q\})$ or to $X \cap (B
\setminus \{q\})$ and obtain a seeing subset of $X$. This is always
possible, unless $q$ fails to see via $X$ some point $a \in X \cap
(A \setminus \{q\})$ and some other point $b \in X \cap (B \setminus
\{q\})$. But then $a$ fails to see $b$ via $X$. (If $[a,b] \subset
X$ then $q \in [a,b]$, as we explained above.) This contradicts our
assumption that $\alpha(X)=2$.
\\

Case 3: $\textrm{dim}K=1$. We will show that in this case $X$ is the union of at most three seeing
subsets. As in case 2, $S$ is the union of two closed convex sets
$A,B$ of full dimension, such that $A \cap B=K$. As $K$ is convex,
$K$ is either a segment,a ray or a line. If $K$ were a line, then
both $A,B$ would be strips or half-planes, and their union would be
convex, which is impossible. So assume $K$ is a segment or a ray.

Suppose $K$ is a ray. W.l.o.g., $K$ lies on the $x$-axis and has a
rightmost point $(w,0)$. If $A$ ($\subset \{(x,y)|y\geq 0\}$) has a
supporting line $L$ at $(w,0)$ that is not horizontal, say
$y=m(x-w)$ $(m \neq 0)$ or $x=w$, then every point $b \in B$ that
lies below the $x$-axis and to the left of $L$ sees via $S$ every
point $a$ of $A$. (The segment $[b,a]$ crosses the $x$-axis within
$K$.) The point $b$ sees, of course, every other point of $B$ via
$S$. Thus $b \in \textrm{ker}S=K$, contrary to our assumption that
$K$ is part of the $x$-axis. By the same token, the set $B$
($\subset \{(x,y)|y\leq 0\}$) does not have a supporting line $L$
through $(w,0)$ that is not horizontal.

At most one of the sets $A,B$ contains points on the $x$-axis to the
right of $(w,0)$. Assume $B$ does not. We claim that $X$ is the
union of two seeing subsets: $B \cap X$ and $(A\setminus K) \cap X$.
We shall first see that $B \cap X$ is a seeing subset of $X$.
Suppose $x,y \in B \cap X$.

If both $x,y$ do not belong to the $x$-axis, take a non-horizontal
line $L$ passing through $(w,0)$ with the points $x,y$ to its right.
As noted before, $L$ does not support $A$, hence there is a point $a
\in A$ to the right of $L$ (see Figure 5). The segment $[x,a]$ meets
the $x$-axis to the right of $(w,0)$. Therefore, $a$ does not see
$x$ via $S$. (Otherwise, $[x,a]$ would be the union of two disjoint
non-empty closed sets $[x,a] \cap A$ and $[x,a] \cap B$.) By the
same token, $a$ does not see $y$ via $S$, and therefore $[x,y]
\subset X$.

If both $x,y$ do belong to the $x$-axis, then $(x,y)
\subset \textrm{rel} \, \textrm{int} K \subset \textrm{int}S \subset X$. If, say, $x \in \textrm{rel} \, \textrm{int} K$ and $y \in B\setminus K$ then
$(x,y) \subset \textrm{rel} \, \textrm{int} B \subset \textrm{int}S \subset X$.

The last case is when, say, $x$ is the endpoint $(w,0)$ of $K$, and
$y \in B\setminus K$. Since $K$ is the only edge of $B$ through
$(w,0)$, $(w,0)$ sees any point $y \in B\setminus K$ via
$\textrm{int}B$, hence via $X$.

We now show that $(A\setminus K) \cap X$ is a seeing subset of $X$.
Take $x,y \in (A\setminus K) \cap X$. Take a non-horizontal line $L$
passing through $(w,0)$ with the points $x,y$ to its right. $L$ does
not support $B$, hence there is a point $b \in B$ to the right of
$L$. According to the considerations brought in the first case
above, $b$ sees via $X$ neither $x$ nor $y$, hence  $[x,y]
\subset X$, so $(A\setminus K) \cap X$ is a seeing subset of $X$ as
well.

Now suppose $K$ is a segment $[u,v]$. If $A$ has a
non-horizontal supporting line at $u$ and a non-horizontal
supporting line at $v$, then there are points in $B\setminus K$ that
are in $\textrm{ker}S$, contrary to our assumption. Therefore, in at
least one of the endpoints of $K$, the only supporting line of $A$
is horizontal. By the same token, in at least one of the endpoints
of $K$, the only supporting line of $B$ is horizontal.

By considerations similar to those brought in the case where
$\textrm{ker}S$ is a ray, the sets $(A\setminus \{(x,y)|y=0\}) \cap
X$, $(B\setminus \{(x,y)|y=0\}) \cap X$ are seeing subsets of $X$.
Now, if $D=\{(x,y)|y=0\} \cap X$ is connected, then we are done as
$X$ is the union of three seeing subsets. Otherwise, $D$ is the
union of two convex sets: $D_1$, the component that includes
$\textrm{rel}\,\textrm{int}K$, and $D_2$. Assume, w.l.o.g., that
$D_2$ is included in $A$ and is disjoint from $B$. Considerations similar to those brought above show that in this case $X$ is the union of two
seeing subsets:
$[(B\setminus \{(x,y)|y=0\}) \cap X] \cup D_1$ and $[(A\setminus \{(x,y)|y=0\})
\cap X] \cup D_2$.
Example 4 shows that in the case where $M_i=\phi$ and $\textrm{dim}K=1$, the number three is best possible:
\\\underline{\textbf{Example 4}}:

Figure 6 describes the set $X$. It is easy to verify that
$\alpha(X)=2$. $\gamma(X)=3$, as there is a 5-circuit of
invisibility.
\\

Case 4: $\textrm{dim}K=2$. This is the most complicated case of the four, which Breen [1974] relates to lengthily.
She claims that in this case, $X$ is the union of four convex sets.
We will show that $X$ is the union of three convex sets. This result
can be viewed as the focal point of the whole paper since it
trivially implies theorem F.
\\\\\underline{Stage 1}: Reduction to the polygonal case:

In this section we intend to show why it is possible to assume that
$S$ is a compact, polygonal set. The result of Lawrence, Hare and
Kenelly [1972] will be useful:

Let $T$ be a subset of a real vector space such that every finite
subset $F\subseteq T$ has a $k$-partition, $\{F_1,...,F_k\}$, with
$\textrm{conv}F_i\subseteq T,\;1\leq i\leq k$. Then $T$ is a union
of $k$ convex sets.

Let $F$ be a finite subset of $X$. We wish to show that $F$ has a $3$-partition, $\{F_1,F_2,F_3\}$,
with $\textrm{conv}F_i\subseteq X$ for $i=1,2,3$. We intend to
construct a set $H$ such that: $\alpha(H)=2$, $\textrm{cl}H$ is
polygonal, $F \subset H \subset X$, $\textrm{cl}H \setminus H
\subset \textrm{bd}\,\textrm{cl}H$ and with
$\textrm{dim\,ker\,cl}H=2$. A representation of $H$ as a union of
three convex sets will imply, in particular, that $F$ has a
partition as required. Therefore, in order to complete our proof, it
will be left to deal with sets $X$ for which $S=\textrm{cl}(X)$ is
polygonal. In the following theorem, we construct a closed set $P$.
We then define $H=P \cap X$ and show that $H$ satisfies the
conditions above.
\begin{theorem}\label{th1}
Suppose $S$ is a closed subset of $\real ^2$, $\alpha(S)\leq 2$,
$\b{0} \in \textrm{int}K$ ($K=\textrm{ker}S$), and $F$ is a given
finite subset of $S$. Then there exists a set $P$ such that:

1) $F \subset P \subset S$

2) $P$ is convex relative to $S$ (hence $\alpha(P) \leq 2$)

3) $\b{0}$ $\in \textrm{int} \; \textrm{ker} P$

4) $P$ is \underline{polygonal}, i.e., $P$consists of a simple
closed polygonal line $\textrm{bd}P$ and its interior.
\end{theorem}

\begin{proof}
We construct the set $P$ in several steps:

 \underline{Step 1}: Add to $F$ the origin $\b{0}$, and, if necessary, a few more points of $S$ (never more than three), so as to make the origin
 $\b{0}$ an interior point of the convex hull of the resulting set. Call the resulting set $F_1$.

 \underline{Step 2}: Define $S_1=S \cap \textrm{conv}F_1$. $S_1$ satisfies all our assumptions on $S$,
and is, in addition, compact. Proceed with $S$ replaced by $S_1$.

 \underline{Step 3}: Replace each point $a \in F_1 \setminus
 \{$\b{0}$\}$ by the intersection of $S_1$ with the closed ray
 $\{\lambda a: \lambda \geq 0 \}$.  Denote the resulting ``sun'' (union of
 segments emanating from $\b{0}$) by $G$. The polygonal set $P$ promised
 in the theorem will be the convex hull of the ``sun'' $G$ relative
 to $S_1$ (or to $S$, doesn't matter).

\underline{Step 4}: Now we start to construct the convex hull of $G$
relative to $S$. Assume $G= \cup_{i=0}^{n-1} [\b{0},a_i]$, where
the points $a_i$ are arranged in order of increasing argument.
Define $a_n=a_0$ and denote by $\Delta_i$ ($i=1,2,...,n$) the
triangle $[\b{0},a_{i-1},a_i]$. For each $i$, $1 \leq i \leq n$, we
define a subset $P_i$ of $\Delta_i$ as follows:

Define:

$\lambda^*=\textrm{max}\{\lambda:0 \leq \lambda \leq 1  \wedge [
\lambda a_{i-1}, a_i] \subset S\}$

 $\mu^*=\textrm{max}\{\mu:0 \leq
\mu \leq 1  \wedge [a_{i-1},\mu a_i] \subset S\}$

The maxima do exist, since $S$ is closed.

Define $P_i=[\b{0},\lambda^* a_{i-1},a_i] \cup [\b{0},a_{i-1},\mu^*
a_i] \subset \Delta_i$. If $[a_{i-1},a_i] \subset S_1$, then
$\lambda^*=\mu^*=1$ and $P_i=\Delta_i$. If not, then $0< \lambda^*
<1$ and $0<\mu^* <1$. ($\lambda^*$ and $\mu^*$ are strictly
positive, since an initial subinterval of $[\b{0},a_{i-1}]$ (and of
$[\b{0},a_i]$) lies in $\textrm{ker}S_1$) In this case , the
intervals $[\lambda^* a_{i-1},a_i]$ and $[a_{i-1}, \mu^* a_i]$ cross
at some point $w_i \in \textrm{int}\Delta_i$, and we obtain:
$P_i=[\b{0}, a_{i-1},w_i] \cup [\b{0},a_i,w_i]$ (see Figure 7).
\begin{claim}\label{cl1}
The set $P_i$ is convex relative to $S$.
\end{claim}
\begin{proof}
This is obvious when $P_i= \Delta_i$. Assume, therefore, that $P_i
\neq \Delta_i$, i.e., $P_i=[\b{0},a_{i-1},w_i] \cup [\b{0},a_i,w_i]=
\Delta_i \setminus (\textrm{int}[a_{i-1},a_i,w_i] \cup
(a_{i-1},a_i))$. Suppose, on the contrary, that some two points $x,y
\in P_i$ see each other via $S$, but not via $P_i$. It follows that,
say, $x \in [\b{0},a_{i-1},w_i]$, $y \in [\b{0},a_i,w_i]$, and the
segment $[x,y]$ passes through $\Delta_i \setminus P_i$ ($=
\textrm{int}[w_i,a_{i-1},a_i] \cup (a_{i-1},a_i)$).

The segment $[x,y]$ cannot meet $(a_{i-1},a_i)$, unless $x=a_{i-1}$ and $y=a_i$, in which case $P_i=\Delta_i$, contrary to our assumption.
 It follows that the segment $[x,y]$ crosses $[w_i,a_{i-1}]$ at some point $x' \neq w_i$ and $[w_i,a_i]$ at some point $y' \neq w_i$.
If $x=x'=a_{i-1}$ and $y' \neq a_i$, extend the segment $[x,y']$
beyond $y'$ into $P_i$, until it hits $[\b{0},a_i]$ at
 some point $y''$ (see the right side of Figure 8). We find that $y''=\mu a_i$ for some $\mu^* < \mu<1$, but $a_{i-1}$ does see $y''$ via $S$,
 contrary to our definition of $\mu^*$.
 We obtain the same type of contradiction when $y=y'=a_i$ but $x' \neq a_{i-1}$.

 Now suppose $y' \neq a_i$, $x' \neq a_{i-1}$. In this case $x' \in (w_i,a_{i-1})$, $y' \in (w_i,a_i)$. Put $z=1/2 (x'+y')$.
 If $a_{i-1}$ sees $z$ via $S$, then it sees via $S$ some point beyond $\mu^* a_i$ on $[\b{0},a_i]$, which is impossible.
 (Note that $S \cap \Delta_i$ is starshaped with respect to $\b{0}$)

 We conclude that $[a_{i-1},z] \nsubseteq S$. By the same token, $[a_i,z] \nsubseteq S$. But $[a_{i-1},a_i]\nsubseteq  S_1$, as well, since
 $P_i \neq \Delta_i$. This contradicts our assumption that $\alpha(S) \leq 2$. (see the left side of Figure 8)
\end{proof}

\underline{Step 5}: Define $P=\cup_{i=1}^n P_i$. Let us check that
$P$ satisfies the requirements of Theorem \ref{th1}.

By our construction, $F \subset F_1 \subset G \subset P \subset S_1
\subset S$.

To prove that $P$ is convex relative to $S$, we take two points $x,y
\in P$ that see each other via $S$, and show that $[x,y] \subset P$.

If $x$ and $y$ belong to the same part $P_i$, then $[x,y] \subset
P_i \subset P$, by Claim \ref{cl1}.

If $\b{0} \in [x,y]$, then $[x,y] \subset P$, since $P$ is
starshaped with respect to $\b{0}$.

Assume, therefore, that $x \in P_i$ and $y \in P_j$, $i < j$, and
that the line through $x,y$ does not pass through the origin. Note
that both $x$ and $y$ lie in $S_1$ ($=S \cap \textrm{conv}F_1$) and
therefore $[x,y] \subset S$ implies $[x,y] \subset S_1$.

For $\nu=0,1,...,n$, denote by $R_{\nu}$ the ray emanating from
$\b{0}$ through $a_{\nu}$ ($R_{\nu}=\{ \lambda a_{\nu} : \lambda
\geq 0\}$). The segment $[x,y]$ crosses the rays $R_i,
R_{i+1},...,R_{j-1}$ (or $R_j, R_{j+1},...,R_n,R_1,...R_{i-1}$) in
this order. Assume, for the sake of of simpler notation, that it
crosses $R_i, R_{i+1},...,R_{j-1}$.

Assume that $[x,y]$ meets $R_{\nu}$ at the point $b_{\nu}=
\lambda_{\nu}a_{\nu}$, where $\lambda_{\nu}>0$. If $\lambda_{\nu}>1$
then $b_{\nu} \notin S_1$, since $a_{\nu}$ is the last point of
$S_1$ on $R_{\nu}$. It follows that $0< \lambda_{\nu} \leq 1$, and
therefore $b_{\nu} \in S_1$, hence $b_{\nu} \in P_{\nu} \cap
P_{\nu+1}$ for $\nu=i, i+1,...,j-1$. Thus $[x,y]=[x,b_i] \cup
[b_i,b_{i+1}] \cup...\cup [b_{j-2},b_{j-1}] \cup [b_{j-1},y]$. By
Claim \ref{cl1}, $[x,b_i] \subset P_i$, $[b_{j-1},y] \subset P_j$
and $[b_{\nu-1},b_{\nu}] \subset P_{\nu}$ for $i<\nu<j$, hence
$[x,y] \subset P$.

To show that $\b{0}$ $\in \textrm{int} \; \textrm{ker} P$, note that
$\b{0}$ $\in \textrm{int} \; \textrm{ker} S$ and $\b{0} \in
\textrm{int}P$. Let $U$ be a neighborhood of $\b{0}$ that lies in $P
\cap \textrm{ker}S$. Every point $u \in U$ sees every point $p \in
P$ ($\subset S$) via $S$, and therefore via $P$, since $P$ is convex
relative to $S$.

Finally, note that the number of edges of the boundary of $P$ never
exceeds $2|F_1|$.
\end{proof}

We can now define the set $H$ as follows: $H = P \cap X$. Let us
show that $H$ satisfies our requirements: (Recall that we need
 $H$ such that: $F \subset H \subset X$, $\alpha(H)=2$, $\textrm{cl}H$ is
polygonal, $\textrm{cl}H \setminus H \subset
\textrm{bd}\,\textrm{cl}H$ and with $\textrm{dim\,ker\,cl}H=2$.)

According to our construction, $F \subset H \subset X$.

Let us show that $H$ is convex relative to $X$: Take two points $a,b
\in H$ such that $[a,b] \subset X$. $a,b \in P$, $[a,b] \subset S$,
so since $P$ is convex relative to $S$, $[a,b] \subset P$. Hence,
$[a,b] \subset X \cap P = H$. Therefore, $\alpha(H) \leq 2$.

$\textrm{int}P \subset \textrm{int}S \subset X$, so $\textrm{int}P
\subset H=P \cap X$. Since $P=\textrm{cl\,int}P$, we find that
$\textrm{cl}H= P$ is polygonal, and $\textrm{cl}H \setminus H
\subset P \setminus \textrm{int}P =\textrm{bd}P =
\textrm{bd}\,\textrm{cl}H$.


Finally, $\b{0}$ $\in \textrm{int} \; \textrm{ker} P$, so
$\textrm{dim\,ker\,cl}H=2$. This concludes the reduction to the
polygonal case. Therefore, we may assume that $S=\textrm{cl}X$ is
polygonal.
\\\\\underline{Stage 2}: Notations:
$S$ is a compact polygonal set. Let $q_1,...,q_n$ be the points of
$Q$ (the lnc points of $S$) ordered in clockwise direction along
$\textrm{bd}(\textrm{conv}Q)$. We assume, for the moment, that
$n\geq 3$. The simpler cases $n=0,1,2$ will be considered
afterwards. $\textrm{conv}Q$ is a polygon with vertices
$q_1,...,q_n$ and edges $e_i=[q_i,q_{i+1}],\;i=1,...,n$ (where
$q_{n+1}=q_1$). By $e_{i_+}$ we denote the closed half-plane
determined by $\textrm{aff} \,e_i$ that misses
$\textrm{int\,conv}Q$. According to Valentine's proof [1957], $S$ is
the union of $\textrm{conv}Q$ and $n$ 'bumps' $W_1,...,W_n$, where
$W_i=S \cap e_{i_+}$. We shall refer to $W_1,...,W_n$ as the
\textbf{leaves} of $S$. Each $W_i$ is a convex polygon and so is the
union $W_i \cup \textrm{conv}Q$, for $i=1,...,n$. Actually, the
union of $\textrm{conv}Q$ with any set of leaves not containing two
adjacent leaves, is a convex polygon.

If we orient the boundary of $S$ clockwise, the boundary of each
leaf $W_i$ (excluding the last edge $e_i$) becomes a directed
polygonal path, with a first edge starting at $q_i$ and a last edge
ending at $q_{i+1}$. Take $l_i$ to be the line spanned by the last
edge of $W_{i-1}$, and $m_i$ to be the line spanned by the first
edge of $W_{i+1}$.
Notice that if $\alpha, \beta, \gamma$ are the angles subtended by
$W_{i-1}$, $\textrm{conv}Q$, and $W_i$ at $q_i$, as in Figure 9,
then the following holds:
\\$\alpha+\beta \leq 180^\circ,\;\;\beta+ \gamma \leq180^\circ$ and $\alpha+ \beta + \gamma>180^\circ$. Therefore, $l_i$ passes either through $\textrm{int}W_i$
or through the basis $e_i$. The same holds for $m_i$. (See Figure
10.)

Denote by $l_{i_+}$ the closed half-plane determined by $l_i$, that
misses $\textrm{conv}Q$ and by $m_{i_+}$ the closed half-plane
determined by $m_i$, that misses $\textrm{conv}Q$.

Done with the description of $S$, we move on to describe $X$: Define
for $i=1,...,n$, $A_i=(W_i \cap X) \setminus Q$. These are the
'leaves' of $X$. (Note that $A_i$ includes the relative interior of
$e_i$, but not its endpoints $q_i,q_{i+1}$.)
 Now, since $\textrm{int}S \subset X$, $X$ can be
represented as the following disjoint union:
\\$X=\textrm{int}(\textrm{conv}Q) \cup (\cup_{i=1}^n{A_i}) \cup (Q
\cap X)$.

*************************************************************
Before entering into more technicalities, we would like to give the
reader an idea of how we properly color $X$ with three colors.

In the original proof of Valentine's Theorem (for $S$), each leaf
$W_i$ (or more precisely, $W_i \setminus e_i$), is colored
uniformly, and two adjacent leaves get different colors. The central
part, $\textrm{conv}Q$, is part of $\textrm{ker}S$, and need not be
colored at all. Thus, two colors suffice locally, and the third
color is only needed to close the circuit when $n$ is odd.

Passing to $X$, the set $A_i$ (=$X \cap W_i \setminus Q$) may miss
some boundary points of $W_i$, and fail to be convex. This
necessitates more than one color for $A_i$. We pass through each
leaf $A_i$ (of $X$) the line $l_i$, that divides $A_i$ into an upper
left part NE ($C_i \cup D_i$ in Figure 10) and a lower right part SW
($F_i \cup E_i$ in Figure 10). The precise definition of this
separation (i.e., which part includes $A_i \cap l_i$) will be given
below.

The NE part is convex, and consists precisely of all points $x \in
A_i$ that fail to see via $X$ some points in $A_{i-1}$. The SW part
is also convex, except (possibly) for some local invisibilities on
the boundary.

We color each of these two parts (NE and SW) uniformly with
different colors. We also have to keep in mind that the color
assigned to NE should be different from the colors assigned to the
adjacent leaf $A_{i-1}$.

Such a coloring will also take care of at least part of the
invisibilities along the boundary of $A_i$. If there is some
invisibility left within the SW part,(this can happen only if the
lines $l_i,m_i$ do not cross within $A_i$), then we fix the coloring
along the boundary using the third color (see the set $G_i$ below).

We can play the same trick with the line $m_i$, coming from the
right, instead of $l_i$, but we shall not use this option. The
following point, however, is important: Whenever the line $l_i$, or
$m_i$, happens to be "horizontal", i.e., coincides with
$\textrm{aff}e_i$, the leaf $A_i$ must be convex, and we may color
it uniformly.

After having colored the $A_i$'s, we finish the job by coloring $Q$.
Each point $q_i \in Q$ belongs to $\textrm{ker}S$, and sees via
$\textrm{int}S$ (hence via $X$) almost all of $X$. $q_i$ may fail to
see via $X$ only points that lie on the last edge of $A_{i-1}$ or
$A_{i-2}$, or on the first edge of $A_{i}$ or $A_{i+1}$ (actually,
on most two of these edges simultaneously). We shall see to it that
the points that $q_i$ does not see use at most two colors, so there
is a third color left for $q_i$.
*************************************************************

Next, we define two partitions of $A_i$ into two parts:\\
\\${l_i}_{(+)}= \begin{cases}
A_i \cap \textrm{int}(l_{i_+})& \text{if $X \cap l_i$ is convex},\\
A_i \cap l_{i_+}& \text{otherwise}.
\end{cases}\;\;\;\;$,     $\;\;\;\;\;l_{i_{(-)}}=A_i \setminus l_{i_{(+)}}$\\
\\\\$m_{i_{(+)}}= \begin{cases}
A_i \cap \textrm{int}(m_{i_+})& \text{if $X \cap m_i$ is convex},\\
A_i \cap m_{i_+}& \text{otherwise}.
\end{cases}\;\;\;$,    $\;\;\;\;\;m_{i_{(-)}}=A_i \setminus m_{i_{(+)}}$ \\
\\\\We shall now see that $l_{i_{(+)}}$ is convex:

If $X \cap l_i$ is not convex then $l_{i_{(+)}}$ includes $l_i \cap A_i$ and there is a point in the last edge of $W_{i-1}$ which is in $X$ and
 does not see any point in $l_{i_{(+)}}$, so due to $\alpha(X)=2$, every two points in $l_{i_{(+)}}$ see each other via $l_{i_{(+)}}$. Otherwise, if
  $X \cap l_i$ is convex then $l_{i_{(+)}}$ does not include $l_i \cap A_i$, so for any
two points $a,b$ in $\;l_{i_{(+)}}$ there is a point in $A_{i-1}$
(close enough to the last edge of $W_{i-1}$) which sees neither $a$
nor $b$, hence $[a,b] \subset l_{i_{(+)}}$.
Similarly, $m_{i_{(+)}}$ is convex as well.

 It is easy to see that any point in $l_{i_{(-)}}$ sees all points in $A_{i-1}$ via $X$.
 (The union $U_i=W_{i-1} \cup \textrm{conv}Q \cup (l_{i_-} \cap W_i)$ is locally convex, and therefore a convex polygon, by Tietze's Theorem. Since
 $\textrm{int}S \subset X$, the only possible invisibilities  in $X \cap U_i$ are along boundary edges of $U_i$. The edge determined by $l_i$
 is taken care of by the exact definition of $l_{i_{(-)}}$. In case $n=3$ there may be another boundary edge of $U_i$ that reaches from
 $A_i$ to $A_{i-1}$, namely, the edge determined by the line $\textrm{aff}(q_{i+1}, q_{i+2})$ (=$\textrm{aff}(q_{i+1},q_{i-1})$).
 If this edge contains a point $x \in A_i$ and a point $y \in A_{i-1}$ then both $x$ and $y$ fail to see via
 $X$ any point $z \in W_{i+1}$, and therefore, $x$ sees $y$ via
 $X$.)
  Similarly, any point in $m_{i_{(-)}}$ sees all points in $A_{i+1}$ via $X$.

Now define:
\\$D_i= l_{i_{(+)}} \cap m_{i_{(+)}}$
\\$C_i= l_{i_{(+)}} \cap m_{i_{(-)}}$
\\$E_i=l_{i_{(-)}} \cap m_{i_{(+)}}$
\\$F_i=l_{i_{(-)}} \cap m_{i_{(-)}} \cap \textrm{ker}X$
\\$G_i=l_{i_{(-)}} \cap m_{i_{(-)}} \setminus \textrm{ker}X$. (see Figure 10)


Notice that Figure 10 describes the case where $l_i,m_i$ meet in
$\textrm{int}A_i$ (then $D_i \neq \phi$ and $G_i=\phi$). $G_i$ may
be non-empty when $l_i, m_i$ do not meet within $A_i$ (as in Figure
11), or even when they meet on the boundary of $A_i$.
\\\\\underline{Stage 3}: The requirements:
As we are going to color each of $C_i \cup D_i,E_i$ uniformly, we
will first see why each of them is convex:

Recall that $D_i= l_{i_{(+)}} \cap m_{i_{(+)}}$, $C_i= l_{i_{(+)}} \cap m_{i_{(-)}}$. So $\;C_i \cup D_i=l_{i_{(+)}}$, which is convex. Now,
$E_i=l_{i_{(-)}} \cap m_{i_{(+)}}\;$ so either $\;E_i=l_{i_-} \cap
m_{i_{(+)}}\;$ or $\;E_i=\textrm{int}(l_{i_-}) \cap m_{i_{(+)}}$. In
any case, $E_i\;$ is convex, as the intersection of two convex sets.

Next, we need to check what are the other requirements for a
coloring $c:X\rightarrow \{0,1,2\}$:
\\\\\textbf{Within each leaf} $A_i:$

$F_i \subset \textrm{ker}(X)$, hence can be given any color. It is
left to check the requirements for its complement in $A_i$.
Since $W_i$ is convex, and $\textrm{int}W_i \subset A_i$,
invisibility within $A_i$ can occur only along edges of $W_i$.
Indeed, two points $a,b \in A_i$ do not see each other via $X$ iff:
\\$(i)$ Both $a$ and $b$ belong to an edge $e$ of $W_i$ (not the
base edge $[q_i,q_{i+1}]$, of course), and
\\$(ii)$ the intersection $e \cap A_i$ is not convex, and
\\$(iii)$ $a,b$ belong to different components of $e \cap A_i$.
\\

If $l_i$ and $m_i$ cross in $\textrm{int}W_i$, then invisibility
within $A_i$ can be along at most one edge $e$ of $W_i$, that goes
all the way from $C_i$ to $E_i$, with one component of $e \cap A_i$
in $C_i$, and the other one in $E_i$. Our coloring will take care of
this invisibility if we require: \\\underline{Requirement 1:} $c(C_i
\cup D_i) \neq c(E_i)$

If $l_i$ and $m_i$ do not cross in $\textrm{int}W_i$, then
invisibility within $A_i$ can occur within edges of $W_i$ that are
not entirely confined to $C_i$ or to $E_i$, i.e., edges that cross
from $C_i$ to $l_{i_(-)} \cap m_{i_(-)}$ or lie entirely in
$l_{i_(-)} \cap m_{i_(-)}$, or cross from  $l_{i_(-)} \cap
m_{i_(-)}$ to $E_i$ (or, possibly, a single edge that reaches  from
$C_i$ through $l_{i_(-)} \cap m_{i_(-)}$ all the way to $E_i$).
$G_i$ consists of all points $a \in A_i$ that belong to $l_{i_(-)}
\cap m_{i_(-)}$ and fail to see some other points $b \in A_i$.

A detailed recipe for a 3-coloring that takes care of all these
invisibilities is given in the next section.
\\\\\textbf{Between two adjacent leaves:}
Two points in adjacent leaves, $a \in A_i,b \in A_{i+1}$, may not
see each other. This can happen only if $a \in m_{i_{(+)}}=D_i \cup
E_i$ and $b \in l_{{i+1}_{(+)}}=C_{i+1} \cup D_{i+1}$. Therefore we
require: \\\underline{Requirement 2:} For each $i, $ $c(E_i) \neq
c(C_{i+1} \cup D_{i+1}) \wedge c(C_i \cup D_i) \neq c(C_{i+1} \cup
D_{i+1})$. (Where $n+1\equiv 1$.)
\\\\\textbf{Involvement of an lnc point $q_i$:}
\\@@@@@@@@@@@@@@@@@@@@@@@@@@@@@@@@@@@@@@@

 $q_i$ may fail to see a point that is in one of the following locations: A point in the last edge
of $A_{i-1}$ ($\subset D_{i-1} \cup E_{i-1}$), a point in the first
edge of $A_i$ ($\subset C_i \cup D_i$), a point in the last edge of
$A_{i-2}$ ($\subset D_{i-2} \cup E_{i-2}$) (this can happen only if
$\,l_{i-1}=\textrm{aff}(q_{i-1},q_i$), in which case $A_{i-1}$ is
convex), or a point in the first edge of $A_{i+1}$ ($\subset C_{i+1}
\cup D_{i+1}$) (this can happen only if
$\,m_{i}=\textrm{aff}(q_{i},q_{i+1}$), in which case $A_i$ is
convex). Now assume that $q_i$ does not see two points $a,b \in X$
at two different locations. Since $\alpha(X)=2$, there are 3 cases
that cannot occur: $a \in A_{i-1},b \in A_i$, $\;a \in A_{i-2},b \in
A_{i-1}$ and $a \in A_{i},b \in A_{i+1}$. This leaves three possible
cases:
\\1. $a \in A_{i-2},b \in A_{i+1}$
\\2. $a \in A_{i-2},b \in A_i$
\\3. $a \in A_{i-1},b \in A_{i+1}$.

Thus it is impossible that $q_i$ won't see points at three different
locations. Stage 4 will show a coloring for $q_i$ that copes with
all three possible cases.
\\\\\textbf{Between two leaves that are not adjacent:}
If a point $x \in A_{i-1}$ does not see a point $y \in A_{i+1}$,
then $x$ is necessarily on the last edge of $W_{i-1}$ and $y$ is on
the first edge of $W_{i+1}$, and these edges lie on the same line,
i.e., $m_i=l_i=\textrm{aff}(q_i,q_{i+1})$. In this case, both $x,y$
do not see any point of $A_i$. This leads to a contradiction to
$\alpha(X)=2$. Therefore, invisibility is impossible among two
leaves of $X$ which are two edges apart.

In any other case, for any two points $a,b \in X$ such that $a
\in A_i, b \in A_{j}$ and $i+2<j<i+n-2$, the segment $(a,b)$ is in
$\textrm{int}S$ (according to Valentine [1957]), and therefore is in
$X$.
\\\\\underline{Stage 4}: An algorithm for a coloring $c:X\rightarrow \{0,1,2\}$:
Since $\textrm{int(conv}Q) \subset \textrm{ker}X$, we only need to
show the coloring of $(\cup_{i=1}^{n}{(A_i\setminus F_i}) \cup(Q
\cap X)$. We start by coloring $A_i\setminus F_i$:
\\\textsc{\emph{\textbf{\underline{\underline{The coloring of $C_i \cup D_i,
E_i$:}}}}}
\\\\\underline{If $n=0$ (mod $3$):} \\$\forall \;1 \leq i \leq n$
$\;\;c(C_i \cup D_i)= i$ $(mod$ $3)$, $c(E_i) = i+2$ $(mod$ $3)$
\\\underline{If $n=1$ (mod $3$):}\\$\forall \;i \neq n$ $\;\;c(C_i \cup D_i)= i$ $(mod$ $3)$, $c(C_n \cup D_n)= 2$. $\forall \;i \neq n-1$
$\;\;c(E_i) = i+2$ $(mod$ $3)$, $c(E_{n-1}) = 1$. See Figure 12 for an example where $n=4$.
\underline{If $n=2$ (mod $3$):}\\$\forall \;1 \leq i \leq n$ $\;\;c(C_i \cup D_i)= i$ $(mod$ $3)$. $\forall \;i \neq n$
$\;\;c(E_i) = i+2$ $(mod$ $3)$, $c(E_{n}) = 0$. See Figure 13 for an example where $n=5$.

We should verify that requirements 1,2 are fulfilled:
\\\textbf{Requirement 1}- For each $i$, $c(C_i \cup D_i) \neq c(E_i)$:
\\\underline{If $n=0$ (mod $3$)}: For $1 \leq i \leq n$, $c(C_i \cup D_i)$= $i$ (mod $3$)$ \neq
c(E_i)=i+2$ $(mod$ $3$)
\\\underline{If $n=1$ (mod $3$)}: For $1 \leq i \leq n-2$, $c(C_i \cup D_i)$= $i$ (mod $3$)$ \neq
c(E_i)=i+2$ $(mod$ $3$). For $i=n-1$, $c(C_{n-1} \cup D_{n-1})= n-1$
$(mod$ $3)=0 \neq c(E_{n-1}) =1$. For $i=n$, $c(C_n \cup D_n)= 2
\neq c(E_n) =n+2$ (mod $3$) $=0$.
\\\underline{If $n=2$ (mod $3$)}: For $i \neq n-1$, $c(C_i \cup D_i)$= $i$ (mod $3$)$ \neq
c(E_i)=i+2$ $(mod$ $3$). For $i=n$, $c(C_n \cup D_n)=n$ (mod $3$)
$=2 \neq c(E_n) =0$.
\\\\\textbf{Requirement 2}- Requirement 2 is composed of two parts:
\\$\bullet$ For each $i$, $c(E_i) \neq c(C_{i+1} \cup D_{i+1})$: (Where $n+1\equiv 1$.)

\underline{If $n=0$ (mod $3$)}: For $1 \leq i \leq n$, $c(E_i) =
i+2$ $(mod$ $3$) $\neq c(C_{i+1} \cup D_{i+1}) =i+1$ (mod $3$)

\underline{If $n=1$ (mod $3$)}:
\\For $1 \leq i \leq n-2$, $c(E_i) =
i+2$ $(mod$ $3$) $\neq c(C_{i+1} \cup D_{i+1}) =i+1$ (mod $3$).
\\For $i=n-1$, $c(E_{n-1}) = 1
\neq c(C_n \cup D_n) =2$
\\For $i=n$, $c(E_n) = n+2$ $(mod$ $3$)$=0
\neq c(C_{1} \cup D_{1}) =1$ (mod $3$)$=1$.

\underline{If $n=2$ (mod $3$)}:
\\For $1 \leq i \leq n-1$, $c(E_i) = i+2$ $(mod$ $3$) $\neq c(C_{i+1}
\cup D_{i+1}) =i+1$ (mod $3$).
\\For $i=n$, $c(E_n) =0
\neq c(C_1 \cup D_1) =1$ (mod $3$)$=1$.
\\\\$\bullet$ For each $i$, $c(C_i \cup D_i) \neq c(C_{i+1}\cup
D_{i+1})$: (Where $n+1\equiv 1$.)

\underline{If $n=0$ (mod $3$)}: For $1 \leq i \leq n$, $c(C_i \cup
D_i)=i$ (mod 3) $\neq c(C_{i+1}\cup D_{i+1})=i+1$ (mod 3).

\underline{If $n=1$ (mod $3$)}:
 \\For $1 \leq i \leq n-2$, $c(C_i \cup
D_i)=i+1$ (mod 3) $\neq c(C_{i+1}\cup D_{i+1})=i+1$ (mod 3).
\\For $i=n-1$, $c(C_{n-1} \cup
D_{n-1})=n-1$ (mod 3) $=0 \neq c(C_n\cup D_n)=2$.
\\For $i=n$, $c(C_n \cup
D_n)=2 \neq c(C_1\cup D_1)=1$.

\underline{If $n=2$ (mod $3$)}:
 \\For $1 \leq i \leq n-1$, $c(C_i \cup
D_i)=i$ (mod 3) $\neq c(C_{i+1}\cup D_{i+1})=i+1$ (mod 3).
\\For $i=n$, $c(C_n \cup
D_n)=n$ (mod 3) $=2 \neq c(C_1\cup D_1)=1$.
\\\\\textsc{\emph{\textbf{\underline{\underline{The coloring of $q_i \in X$:}}}}}

As mentioned in stage 3, there are three possible 'maximal' cases where $q_i$
does not see via $X$ two points $\{a,b\}$ at two different
locations:
 \\\textbf{Case 1}: $a \in A_{i-2},\;b \in A_{i+1}$: In
this case, $a \in D_{i-2} \cup E_{i-2}$ and $b \in C_{i+1} \cup
D_{i+1}$. In addition, $a \in \textrm{aff}(q_{i-1},q_{i})$, $b \in
\textrm{aff}(q_{i},q_{i+1})$. It follows that there are no two
adjacent lnc points $q_i, q_{i+1}$ such that both points satisfy case 1.
Otherwise, if $q_i$ does not see via $X$ two points $a \in
A_{i-2},\;b \in A_{i+1}$ and $q_{i+1}$ does not see via $X$ two
points $a' \in A_{i-1},\;b' \in A_{i+2}$ , then $a' \in
\textrm{aff}(q_i,q_{i+1})$, $b \in \textrm{aff}(q_{i},q_{i+1})$.
Hence, the points $a',q_i,b$ do not see each other via $X$, a
contradiction to $\alpha(X)=2$. We shall take care of case 1 with regards to the three cases of $n$ mod $3$:

\underline{If $n=0$ (mod $3$)}: According to the coloring $c$, in
this case, for every $1 \leq i \leq n$, $\;\;c(C_{i+1} \cup
D_{i+1})=i+1 (mod \;\;3)$, $c(C_{i-2} \cup D_{i-2})=i-2 (mod\;\;3)$,
therefore $c(C_{i+1} \cup D_{i+1})=c(C_{i-2} \cup D_{i-2})$ so those
two sets are colored by the same color. The set $E_{i-2}$ is colored
by a second color. As these sets are the only constraints for $q_i$,
we will color $q_i$ by the third color that is left 'free'.

\underline{If $n=1$ (mod $3$)}:

For $3 \leq i \leq n-2$, the same considerations as in the case
$n=0$ (mod $3$) hold.

For $i=n$, $c(C_{i+1} \cup D_{i+1})=c(C_1 \cup D_1)=1$, $c(C_{i-2}
\cup D_{i-2})=c(C_{n-2} \cup D_{n-2})=n-2(mod\;\;3)=2$,
$c(E_{i-2})=c(E_{n-2})=n$ (mod $3$)$=1$, so we define $c(q_n)=0$.

For $i=2$, $c(C_{i+1} \cup D_{i+1})=c(C_3 \cup D_3)=0$, $c(C_{i-2}
\cup D_{i-2})=c(C_n \cup D_n)=2$, $c(E_{i-2})=c(E_{n})=n+2$ (mod
$3$)$=0$, so we define $c(q_2)=1$.

For $i=1$ or $i=n-1$, there is no color left to assign to $q_i$,
since all three colors are occupied by $C_{i+1} \cup D_{i+1}$ ,$C_{i-2}
\cup D_{i-2}$ and $E_{i-2}$. We solve this situation by renumbering the
lnc points in a way that case 1 occur neither in $q_1$, nor in $q_{n-1}$.
This is always possible because for every lnc point satisfying case
1, both of the adjacent lnc points do not satisfy case 1 so there
are always to find two lnc points which can be marked as $q_1$ and
$q_{n-1}$.

\underline{If $n=2$ (mod $3$)}:

For $3 \leq i \leq n-1$, the same considerations as in the case
$n=0$ (mod $3$) hold.

For $i=2$, $c(C_{i+1} \cup D_{i+1})=c(C_3 \cup D_3)=0$, $c(C_{i-2}
\cup D_{i-2})=c(C_n \cup D_n)=n$ (mod $3$)$=2$,
$c(E_{i-2})=c(E_{n})=0$, so we define $c(q_2)=1$.

For $i=1$ or $i=n$, there is no color left to assign to $q_i$, since
all three colors are occupied by $C_{i+1} \cup D_{i+1}$, $C_{i-2} \cup
D_{i-2}$ and $E_{i-2}$. If there are to find two adjacent lnc points not satisfying case 1, then  we renumber the
lnc points in a way that case 1 occur neither in $q_1$, nor in $q_n$. Notice that if $n$ is odd, then this is always possible, since for every lnc point satisfying case 1, both of the adjacent lnc points do not satisfy case 1. If $n$ is even, then it might be possible that there are no two adjacent lnc points not satisfying case 1. This can happen only if the lnc points ${q_i}$ satisfy case 1 alternately, and it follows, that in this case all $A_i$ are convex. In this special case we perform a different coloring: we color all $A_i$ alternately by two colors and the lnc points by the third color. (See Figure 14)
\\\textbf{Case 2}: $a \in A_{i-2},\;b \in A_i$: Color $q_i$ by $c(C_{i-1} \cup D_{i-1})$. Notice that $b \in C_i \cup D_i$.
 Now, $c(C_{i-1} \cup D_{i-1}) \neq c(C_i \cup D_i)$, so $q_i$ and $b$ are colored differently. Similarly,
$c(C_{i-1} \cup D_{i-1}) \neq c(E_{i-2})$ and $c(C_{i-1} \cup
D_{i-1}) \neq c(C_{i-2} \cup D_{i-2})$ so $q_i$ and $a$ are colored
differently.
\\\textbf{Case 3}: $a \in A_{i-1},\;b \in A_{i+1}$: Color $q_i$ by $c(C_i \cup D_i)$. In this case by similar considerations to those in case 2, $q_i$
is colored differently from both $a,b$ .
\\\textbf{The coloring of $G_i$:}
$G_i$ is the disjoint union of finitely many connected components
$b_j$. Each component is either a single point, or a line segment,
or the union of two line segments with a common endpoint (This
common endpoint must, of course, be a vertex of $W_i$). Each edge of
$W_i$ meets at most two components of $G_i$. If $G_i \neq \phi$,
then $D_i = \phi$, and all the components $b_j$ lie in the gap
between $C_i$ and $E_i$ (see Figure 11). Number the components
$b_1,...,b_t$ in the order they appear on the boundary of $A_i$,
with $b_1$ closest to $C_i$ and $b_t$ closest to $E_i$.

Points of $b_1$ may fail to see points of $C_i$, and points of
$b_t$ may fail to see points of $E_i$. Beyond that, points of $G_i$
may fail to see each other via $X$ only if they belong to adjacent
components $b_j,\,b_{j+1}$.

Assume $\{0,1,2\}=\{p,q,r\}$, where $p=c(E_i)$ and $q=c(C_i)$. Color
the components $b_j$ as follows:
\\$c(b_j)= \begin{cases}
p& \text{if $j$ is odd and $j \neq t$},\\
q& \text{if $j$ is even},\\
r& \text{if $j$ is odd and $j=t$}.
\end{cases}$

Note that $c(b_1) \neq c(C_i)$, $c(b_t) \neq c(E_i)$ and $c(b_j)
\neq c(b_{j+1})$ for $j=1,2,...,t-1$.
This finishes the description of a 3-coloring of $X$.

It is left to deal with the cases $n=0,1,2$.
\\\\\underline{For $n=0$}: Let us show that in this case, when $M_i= \phi$ and $S$ is convex, then
$\gamma(X) \leq 3$. Take $a,b \in X$. If $[a,b] \cap \textrm{int}S
\neq \phi$ then $[a,b] \subset \textrm{int}S$ and since
$\textrm{int}S \subset X$, $[a,b] \subset X$. Hence, $\textrm{int}S
\subset \textrm{ker}X$. Therefore, we only need to show a coloring
of $X \cap \textrm{bd}S$ with three colors. Since $S$ is polygonal,
$\textrm{bd}S$ is a circuit of edges. Each edge of $S$ meets at most
two components of $X \cap \textrm{bd}S$, so $X \cap \textrm{bd}S$ is
the disjoint union of a finite number of components $e_j$. Each
component is either a single point, or a line segment, or the union
of two line segments with a common endpoint (This common endpoint
must, of course, be a vertex of $S$). Number the components
$e_1,...,e_t$ in the order they appear on $\textrm{bd}S$, where
$e_1$ is adjacent to $e_t$ .

Points of $\textrm{bd}S \cap X$
may fail to see each other via $X$ only if they belong to adjacent
components $e_j,\,e_{j+1}$. Color the components $e_j$ as follows:
\\$c(e_j)= \begin{cases}
1& \text{if $j$ is odd and $j \neq t$},\\
2& \text{if $j$ is even},\\
3& \text{if $j$ is odd and $j=t$}.
\end{cases}$

Note that $c(e_j)
\neq c(e_{j+1})$ for $j=1,2,...,t-1$ and that $c(e_t)
\neq c(e_1)$. Let us mention that in this case (where $\textrm{cl}X=S$ is convex and $M_i=\phi$), $\gamma(X)=3$ iff $S$ is an odd-sided
convex polygon (not a triangle), $X$ contains all vertices of $S$ and misses at least one point in each edge of $S$.
\\\\\underline{For $n=1$}: $Q=\{q\}$ and $\textrm{bd}S$ is a polygon $<q,
a_1,a_2,...,a_k,q>$.
(see Figure 15). Define $l$ to be the line spanned by $[a_k,q]$ and
$m$ to be the line spanned by $[q,a_1]$. $l_+$ is the closed
half-plane determined by $l$ that contains $a_1$, and $m_+$ is the
closed half-plane determined by $m$ that contains $a_k$. Note that
$\textrm{ker}S=S \cap l_- \cap m_-$. Define:
\\\\$l_{(+)}= \begin{cases}
X \cap \textrm{int}(l_+)& \text{if $X \cap l$ is convex},\\
(X\setminus [a_k,q]) \cap l_+& \text{otherwise}.
\end{cases}$ \\\\and similarly define:
\\\\$m_{(+)}= \begin{cases}
X \cap \textrm{int}(m_+)& \text{if $X \cap m$ is convex},\\
(X\setminus [a_1,q]) \cap m_+& \text{otherwise}.
\end{cases}$

By considerations similar to those appearing in the case $n\geq
3$ (regarding $l_{i_{(+)}}$ and $m_{i_{(+)}}$), both $l_{(+)}$ and
$m_{(+)}$ are convex.

Define $F=X\setminus \{q\} \setminus l_{(+)}\setminus
m_{(+)}$. A glance at Figure 15 shows that $\textrm{int}F \subset
\textrm{ker}X$. A close look at the definition of $l_{(+)}$ and
$m_{(+)}$ shows that every point $x \in F \cap \textrm{int}S$
belongs to $\textrm{ker}X$. There may be some invisibilities in $X$
along edges of $S$ that meet $F$. Define $G=F\setminus
\textrm{ker}X$. $G$ is composed of connected components along $F
\cap \textrm{bd}S$.

We define now the coloring $c$ of $X$: $X \cap \textrm{ker}X$ can be given any color, $c(l_{(+)})=1$,
$c(m_{(+)})=2$. $G$ will be colored in the same manner as $G_i$ is
colored in the case $n \geq 3$ (here three colors might be needed),
and it is left to color $q$: If there is no point in $l_{(+)}$ which
is invisible from $q$ then $c(q)=1$. Otherwise, there is $x \in
[a_1,q]$ that does not see $q$ via $X$. Now, If there is no point in
$m_{(+)}$ which is invisible from $q$ then $c(q)=2$. Otherwise,
there is a point $y \in [a_k,q]$ that doest not see $q$ via $X$. The
points $x,y,q$ form a non-seeing subset of $X$, a contradiction to
$\alpha(X)=2$.
\\\\\underline{For $n=2$}: A similar construction and coloring as in the case $n\geq 3$ is available
here. $Q=\{q_1,q_2\}$, $e=[q_1,q_2]$. The line $\textrm{aff}\,e$
divides the plane into two closed half-planes: $e_+,e_-$ (in Figure
16, we take $e_+$ to be above $e_-$). Denote $W_1=S \cap e_+$ and
$W_2=S \cap e_-$. For $i=1,2$, $\textrm{bd}W_i$ includes, besides
$e$, a sequence of edges $e_{i1},e_{i2},...,e_{i{n_i}}$, ordered in
clockwise direction. Notice that the first edge $e_{i1}$, and the
last edge, $e_{i{n_i}}$, may lie on the same line as $e$. $l_1$ is
the line spanned by the last edge of $\textrm{bd}W_2$, $e_{2{n_2}}$,
and $l_2$ is the line spanned by the last edge of $\textrm{bd}W_1$,
$e_{1{n_1}}$. Take  $m_1$ to be the line spanned by $e_{21}$, and
$m_2$ to be the line spanned by $e_{11}$. (See Figure 16.) Denote by
$l_{1_+}$ the closed half-plane determined by $l_1$, not containing
$W_2$, and by $l_{2_+}$ the closed half-plane determined by $l_2$,
not containing $W_1$. Similarly, $m_{1_+}$ is the closed half-plane
determined by $m_1$, not containing $W_2$, and $m_{2_+}$ is the
closed half-plane determined by $m_2$, not containing $W_1$. Done
with the description of $S$, we move on to describe $X$. For
$i=1,2$, $A_i=(W_i \cap X)\setminus e$, $X=A_1 \cup A_2 \cup
(q_1,q_2) \cup (Q \cap X)$. Each $A_i$ has the following partitions:
\\\\${l_i}_{(+)}= \begin{cases}
A_i \cap \textrm{int}(l_{i_+})& \text{if $X \cap l_i$ is convex},\\
A_i \cap l_{i_+}& \text{otherwise}.
\end{cases}\;\;\;\;$,     $\;\;\;\;\;l_{i_{(-)}}=A_i \setminus l_{i_{(+)}}$\\
\\\\$m_{i_{(+)}}= \begin{cases}
A_i \cap \textrm{int}(m_{i_+})& \text{if $X \cap m_i$ is convex},\\
A_i \cap m_{i_+}& \text{otherwise}.
\end{cases}\;\;\;$,    $\;\;\;\;\;m_{i_{(-)}}=A_i \setminus m_{i_{(+)}}$ \\
\\Now, define:
\\$D_i= l_{i_{(+)}} \cap m_{i_{(+)}}$
\\$C_i= l_{i_{(+)}} \cap m_{i_{(-)}}$
\\$E_i=l_{i_{(-)}} \cap m_{i_{(+)}}$
\\$F_i=l_{i_{(-)}} \cap m_{i_{(-)}} \cap \textrm{ker}X$
\\$G_i=l_{i_{(-)}} \cap m_{i_{(-)}} \setminus \textrm{ker}X$. (see Figure 16)

As in the case $n \geq 3$, ${l_i}_{(+)}$, ${m_i}_{(+)}$ are
convex,
so again, this entails the convexity of $C_i \cup D_i$ and $E_i$ for
$i=1,2$. All points in ${l_1}_{(-)}$ see all points in
${l_2}_{(-)}$. In addition, all points in ${m_1}_{(-)}$ see all
points in ${m_2}_{(-)}$. Therefore, $C_1,C_2$ see each other and
$E_1,E_2$ see each other. All the above enables the same
 coloring as the one described in the case $n \geq 3$, which leads to the
following: $c(C_1 \cup D_1)=1,c(E_1)=3, c(C_2 \cup D_2)=2,c(E_2)=3$.
The sets $G_1,G_2$ will also be colored as in the case $n\geq 3$.
The considerations written above together with those appearing in
the case $n\geq 3$ clarify why this coloring is proper.
We can now conclude that if $M_i=\phi$ and $\textrm{dim}K=2$, then
$\gamma(X)\leq 3$.

Example 5 (due to Breen [1974]) shows that the number three is best possible.
\\\underline{\textbf{Example 5}}:

Let $P$ be a regular pentagon. Define $X=(P \setminus \textrm{bd}P) \cup \textrm{vert}P$ (see Figure 17).

$\alpha(X)=2$: The only points in $X$ that are not in
$\textrm{ker}X$ are the vertices of $P$. The only points that a
vertex does not see via $X$ are the two adjacent vertices, but these
two see each other via $X$.

There is a 5-circuit of invisibility, therefore $\gamma(X)\geq 3$.
\section{Proof of Theorem F}
Assume $M_i=\{(0,0)\}$. Define $A=\{(x,y) \in \real^2:y>0 \vee (y=0
\wedge x>0) \}$, $B=\{(x,y) \in \real^2:y<0 \vee (y=0 \wedge x<0)
\}$. $A \cup B= \real^2\setminus \{(0,0)\}$ and $(0,0) \notin X$,
therefore $X=(X \cap A) \cup (X \cap B)$. $A,B$ are both convex, so
each of the sets $X \cap A, X \cap B$ satisfies the conditions of
theorem E, and therefore each is a union of three convex sets.
Hence, $X$ is the union of six convex sets.

Example 6 shows that the number six is best possible.
\\\underline{\textbf{Example 6}}:

We describe a set $X \subset \real^2$ and show that $\alpha(X)=2$ and $\gamma(X)>5$.
Let $P$ be a regular $48$-gon with center $O$, vertices $p_0,p_1,...,p_{48}$ ($p_0=p_{48}$) and edges $e_i=[p_{i-1},p_i]$ ($i=1,2,...,48$). Above each edge $e_i$ erect a triangular dome $T_i=[q_{i-1},q_i,,t_i]$. The interior angles of each $T_i$ are as follows:
\\At the odd-numbered base vertex: $7.5^\circ$ ($=360^\circ/48$).
\\At the even-numbered base vertex: $6^\circ$.
\\At the tip $t_i$: $166.5^\circ$.

Define $S=P \cup (\cup_{i=1}^{48} T_i)$. $X$ is obtained from $S$  by removing the odd-numbered vertices $p_{2k-1}$ ($k=1,2,...24$) and the center $O$.

Note that each odd-numbered vertex $p_{2k-1}$ is the crossing point of the segments $[p_{2k-2},t_{2k}]$ and $[t_{2k-1},p_{2k}]$.
Moreover, the sum of the interior angles of $T_{2k}$, $P$ and
$T_{2k+1}$ at the even-numbered vertex $p_{2k}$  is
$184.5^\circ(>180^\circ)$. Therefore, $t_{2k}$ and $t_{2k+1}$ do not
see each other via $X$. Thus we see that, for each $k$, the points
$<p_{2k},t_{2k-1},t_{2k},t_{2k+1},t_{2k+2},p_{2k}>$ form a 5-circuit
of invisibility in $X$.
Note also that $t_i$ and $t_{i+2}$ always see each other via $\textrm{int}X$. (See Figure 18)

Next, we show that $\alpha(X)=2$. Note that $S=\textrm{cl}X$ is the
union of the (closed) 48-gon $P$ and 48 triangular (closed) domes,
with $\textrm{ker}S=P$. The union of $P$ and any collection of
non-adjacent (closed) domes is convex. Thus two points of $S$ fail
to see each other via $S$ only if they belong to adjacent domes. A
close look at Figure 18 shows that if $z,w \in X$, and the open
segment $(z,w)$ passes through one of the removed
 vertices of $P$, then $z$ and $w$ must belong to two adjacent domes.

Assume that three points $a,b,c$ form a 3-circuit of invisibility in $X$.
\\a) If $a \in \textrm{int}P$ then all points of $X$ that $a$ does not see via $X$, in particular $b$ and $c$, lie on the ray
 $R_{-a}= \{(1+ \lambda)O-\lambda a: \lambda>0 \}$. The intersection of this ray with $X$ is convex, and thus $[b,c] \subset X$, contrary to our assumption.
 \\Thus we may assume  that each of $a,b$ and $c$ belongs to one of the 48 closed triangular domes $T_1,...,T_{48}$.
 \\b) If $a$ is an even-numbered vertex of $P$, say $p_2$ (see Figure 18), then the only points of $X\setminus \textrm{int}P$ that do not see $a$ via $X$ are: the opposite vertex ($p_{26}$), the points of the segment $[t_1,p_1)$ and the points of the segment $(p_3,t_4]$.
 But all these points see each other via $X$, so again, $[b,c] \subset X$.
 \\c) Assume, therefore, that $a,b,c, \in X\setminus (\textrm{int}P \cup \textrm{vert}P)$. It follows that each of these three points belongs to a
 unique dome (the only points of $X$ that are common to two domes are the even-numbered vertices of $P$). Denote by $T_a,T_b,T_c$ the domes that
 contain $a,b$ and $c$ respectively. If $[a,b] \nsubseteq X$, then $T_a$ and $T_b$ must be adjacent or (if $O\in (a,b)$) opposite domes.
 Same for $T_a$ and $T_c$, $T_b$ and $T_c$. But there are no three different domes such that each two are either adjacent or opposite. Thus $\alpha(X)=2$.

It is left to convince the reader that $\gamma(X)>5$: The 12 rays
$\overrightarrow{OP_{4k}}$
\\($k=0,...,11$) divide $X$ into twelve
congruent sectors. Each sector includes four consecutive edges of
$P$, and the corresponding domes, and contains a 5-circuit of
invisibility. Figure 18 represents one sector. (The central angles
are, of course, exaggerated)

It follows that each sector is not a union of two convex sets, and therefore, in any covering of $X$ by convex subsets, each sector will meet at least three of the covering subsets.
Now, assume to the contrary that $X$ is the union of $5$ convex
sets. Let's try to evaluate the number of incidences between the
five convex sets and the 12 sectors. On the one hand, as every
sector meets at least three convex sets, this number is no less than
$3\cdot 12=36$. On the other hand, as none of the convex sets
includes the center $O$, each convex set lies on one side of a line
through $O$, and therefore meets at the most 7 sectors. Therefore,
the number of incidences is not more than $7\cdot 5=35$, a
contradiction.

\section{Proof of Theorem G}
In the proof of Theorem E, in the case where $n=0$, we gave the
following characterization: When $S$ is convex and $M_i=\phi$,
$\gamma(X)=3$ iff $S$ is an odd-sided convex polygon (not a
triangle), $X$ contains all vertices of $S$ and misses at least one
point in each edge of $S$. (otherwise $\gamma(X)\leq 2$)

Assume $M_i=\{(0,0)\}$. Define $A=\{(x,y) \in \real^2:y>0 \vee (y=0
\wedge x>0) \}$, $B=\{(x,y) \in \real^2:y<0 \vee (y=0 \wedge x<0)
\}$. $A \cup B= \real^2\setminus \{(0,0)\}$ and $(0,0) \notin X$,
therefore $X=(X \cap A) \cup (X \cap B)$. $A,B$ and $\textrm{cl}X$
are all convex.

We show now that the set $\textrm{cl}(X \cap A)$ is convex: Let $a,b
\in \textrm{cl}(X \cap A)$ and assume there is $c \in (a,b)$ such
that $c \notin\textrm{cl}(X \cap A)$. There are two series
$(a_n),(b_n) \subset X \cap A$ such that $(a_n)\rightarrow
a,(b_n)\rightarrow b$. If  $c \notin\textrm{cl}(X \cap A)$ then
there is a neighborhood $U$ of $c$ such that $U \cap (X \cap A)=
\phi$. It follows that there exists $c' \in U$ satisfying the
following: 1. $c' \in (a_n,b_n)$ for some $n$. 2. There is a
neighborhood $U'$ of $c'$ satisfying $ U' \subset U$ and $U' \subset
A$. It follows that $U' \cap X = \phi$, hence, $c' \notin
\textrm{cl} X$, contradiction to the convexity of $\textrm{cl}X$.
Similarly, $\textrm{cl}(X \cap B)$ is convex as well. Therefore,
each of $X \cap A$, $X \cap B$ satisfies the conditions of the
characterization brought above. Each of these sets has an edge with
a missing vertex, therefore, according to the characterization, each
of $X \cap A$, $X \cap B$ is the union of at most two convex sets,
hence $\gamma(X)\leq 4$.

Example 7 shows that the number four is best possible.
\\\underline{\textbf{Example 7}}:

We describe a convex set $X \subset \real^2$ and show that $\alpha(X)=2$ and $\gamma(X)\geq 4$.
Let $P$ be a regular $7$-gon with center $O$. We define
$X=P \setminus (\{O\} \cup \textrm{bd}P) \cup \textrm{vert}P$
($X$ is obtained from $P$ by removing the center $O$ and the
relative interiors of all edges). Let $C$ be a convex subset of $X$.
Since $O \notin X$, $C$ is included in a closed half-plane $H$ with
$O \in \textrm{bd}H$. $H$ intersects $\textrm{vert}P$ in a stretch
of three or four consecutive vertices. But two adjacent vertices of
$X$ do not see each other via $X$. Therefore $C$ contains at most
two vertices of $P$. It follows that $\gamma(X) \geq 4$.

 We show now that $\alpha(X)=2$: $X \setminus \textrm{vert}P$ is
the union of two convex sets. Therefore. if there is a 3-circuit of
invisibility in $X$, then it must contain a vertex of $P$. For each
vertex of $P$, the points in $X$ which it does not see via $X$ are
the two adjacent vertices in $P$ and the opposite radius. Notice
that all these points see each other via $X$. Therefore, a vertex of
$P$ cannot be a part of a 3-circuit of invisibility, hence
$\alpha(X)=2$.\\

It is left to show that if $M_b=\phi$ or $M_b=\textrm{bd}S$, then
$\gamma(X)=2$. Indeed, if $M_b=\phi$ then $X=S\setminus\{O\})=(S\cap
A) \cup (S \cap B)$, where $A,B$ are the convex sets defined above.
If $M_b=\textrm{bd}S$, then $X=\textrm{int}S\setminus
\{O\}=(\textrm{int}S \cap A) \cup (\textrm{int}S \cap B)$.

\section{Proof of Main Theorem 2}
Example 8 shows that the number four is best possible.
\\\underline{\textbf{Example 8}}:

We describe a set $X \subset \real^2$ and show that $\alpha(X)=2$ and $\beta(X)\geq 4$.
\\Let $P$ be a square with center $O$.

$\alpha(X)=2$: The proof is similar to the one in Example 6.

$\beta(X)\geq 4$: Assume to the contrary that $X$ can be colored by
three colors, in a way that for each color, the set of points in $X$
colored by that color, is a seeing subset of $X$. The points
$<1,2,3,4,5>$ form a 5-circuit of invisibility in $X$, therefore,
three colors are needed in order to color them. It follows that in
order to color the set of points $\{a,b,c,d,e\}$, at least two
colors are needed, as otherwise, if we color all of them by one
color, then there must be a point among $\{1,2,3,4,5\}$ colored by
the same color and this point does not see one of the points
$\{a,b,c,d,e\}$.

Similarly, at least two colors are needed in
order to color the points $\{a',b',c',d',e'\}$.
All points in the set $\{a,b,c,d,e\}$ do not see any point in
$\{a',b',c',d',e'\}$, therefore, at least four colors are needed to
color the set $\{a,b,c,d,e,a',b',c',d',e'\}$, a contradiction.

\section{Proof of Main Theorem 3}
If $\beta(X)=2$ then $\beta(S)\leq 2$ so $\alpha(S) \leq 2$. If
$M_i=\phi$ then according to Theorem E, $\gamma(X) \leq 3$. If
$|M_i|>1$ then according to Theorem C, $\gamma(X) \leq 3$.  It is
left to handle the case $|M_i|=1$:

Assume $M_i=\{p\}$. According to Lemma \ref{noa3}, $p \in
\textrm{ker}S$. We wish to show that for every $x \in X$, $(p,x]
\subset X$. Assume to the contrary that there is $x \in X$ such that
$(p,x] \nsubseteq X$. Hence, there is $y \in (p,x]$, $y \notin X$.
$p \in \textrm{int}S$ is the only point in $M_i$, so there is some
neighborhood $U$ of $p$ such that $U\setminus \{p\} \subset X$.
Hence, there are two points $a,b \in \textrm{aff}(p,x) \cap U$ such
that $a,b$ do not see each other via $X$. Due to the absence of the
point $y$ in $X$, both $a,b$ do not see the point $x$ via $X$, and
we have a contradiction to $\alpha(S) \leq 2$.

 $X$ is the union of two seeing subsets $A,B$. We shall see now that $X \cup \{p\}$ is the union of two convex
sets: $\textrm{conv}(A \cup \{p\}),\textrm{conv}(B \cup \{p\})$. In
order to show that, we need to show that $\textrm{conv}(A \cup
\{p\}) \subset X \cup \{p\}$.

We shall use the 'strong' version of Caratheodory: If $C \subset
\real ^d$, $p \in C$ then each point in $\textrm{conv}C$ can be
represented as a convex combination of an affine independent subset
of $C$ containing $p$.

Hence, every point in $\textrm{conv}(A \cup \{p\})$ is a convex
combination of $p$ and other two points $a,b \in A$. $[a,b] \subset
X$ and for every $y \in [a,b]$, $(p,y] \subset X$. Therefore,
$[p,a,b] \subset X \cup \{p\}$. Symmetrically, $\textrm{conv}(B \cup
\{p\}) \subset X \cup \{p\}$, which implies that $X \cup \{p\}$ is
the union of two convex sets. This implies that $X$ is the union of
four convex sets.
\bibliographystyle{amsplain}
\bibliography{}
Breen, M. (1974), ``Decomposition theorems for 3-convex subsets of
the plane'', Pacific J. Math., Vol. 53, 43-57.
\newline \newline Breen, M. and D. Kay (1976), ``General decomposition theorems for m-convex sets
in the plane'', Israel J. Math., Vol. 24, 217-233.
\newline\newline Lawrence J.F., W.R. Hare and J.W. Kenelly
(1972), ``Finite unions of convex sets'' , Proc. Amer. Math. Soc.,
34, 225-228.
\newline\newline Matou$\check{\text{s}}$ek J. and P. Valtr (1999), ``On visibility and covering by convex sets'', Israel J.
Math., Vol. 113, 341-379.
\newline\newline Perles, M.A. and S. Shelah (1990), ``A closed (n+1)-convex set in $\real^2$ is a union of $n^6$ convex sets'', Israel J.
 Math., Vol. 70, 305-312.
\newline\newline Valentine, F.A. (1957), ``A three point convexity property'', Pacific J. Math.,
7(2), 1227-1235.
\newline
\newline

\end{document}